\documentclass[10pt]{article}
\usepackage[dvips]{graphicx}
\usepackage{times}
\usepackage{type1cm}
\usepackage{amsmath,amsthm,amsfonts,amssymb,amscd}

\makeatletter
\renewcommand{\section}{\@startsection{section}{0}{0mm}{10pt}{5pt}{\normalsize \bf }}
\makeatother

\newtheorem{thm}{\indent Theorem}[section]
\newtheorem{cor}[thm]{\indent Corollary}
\newtheorem{defn}[thm]{\indent Definition}

\newtheorem{lem}[thm]{\indent Lemma}
\newtheorem{prop}[thm]{\indent Proposition}

\newtheorem{conj}[thm]{\indent Conjecture}
\newtheorem{rem}[thm]{\indent Remark}
\newtheorem{rems}[thm]{\indent Remarks}

\title{\bf \normalsize COMPUTATIONS OF TURAEV-VIRO-OCNEANU INVARIANTS OF 3-MANIFOLDS 
FROM SUBFACTORS}

\author{
{\footnotesize NOBUYA SATO}
\footnote{Supported in part by the
Grants-in-Aid for Scientific Research, JSPS.}\\
{\footnotesize \it Department of Mathematics and Information Sciences}\\
{\footnotesize \it Osaka Prefecture University, Sakai, Osaka, 599-8531, JAPAN}\\
{\footnotesize \it e-mail:nobuya@mi.cias.osakafu-u.ac.jp}\\ 
 \\ 
{\footnotesize MICHIHISA WAKUI}\\
{\footnotesize \it Department of Mathematics, Osaka University, Toyonaka, Osaka, 560-0043, JAPAN}\\
{\footnotesize \it e-mail:wakui@math.sci.osaka-u.ac.jp}
\\ 
}

\begin{document}

\maketitle
\par \ 
\par 
\centerline{\normalsize ABSTRACT}
\par \ 

\begin{minipage}{4.5in}
\footnotesize 
\baselineskip=10pt  
\qquad  
In this paper, we establish a rigorous correspondence between the two tube algebras, 
that one comes from the Turaev-Viro-Ocneanu TQFT introduced by Ocneanu 
and another comes from the sector theory introduced by Izumi, 
and construct a canonical isomorphism between the centers 
of the two tube algebras, which is a conjugate linear isomorphism 
preserving the products of the two algebras 
and commuting with the actions of $SL(2,\mathbb{Z})$. 
Via this correspondence and the Dehn surgery formula, 
we compute Turaev-Viro-Ocneanu invariants 
from several subfactors for basic $3$-manifolds including 
lens spaces and Brieskorn $3$-manifolds by using Izumi's data written in terms of sectors.
\end{minipage} 

\par \ 
\baselineskip=13pt 

\section{\kern-1em .\ \ Introduction} 
\par 
At the beginning of the 1990's, a $(2+1)$-dimensional 
unitary topological quantum field theory, in short, TQFT, was 
introduced by A. Ocneanu \cite{Ocneanu1} by using a type II$_1$ subfactor with finite index 
and finite depth as a generalization of the Turaev-Viro TQFT \cite{TV} 
which was derived from the quantum group $U_q(sl(2,\Bbb{C}))$ at certain roots 
of unity. 
We call such a TQFT a Turaev-Viro-Ocneanu TQFT. 
K. Suzuki and the second author \cite{SuzukiWakui} have found a 
Verlinde basis in the sense of \cite{KSW} for the Turaev-Viro-Ocneanu TQFT from an $E_6$-subfactor, and 
computed the invariant for basic $3$-manifolds including lens spaces 
$L(p,q)$, where $p, q$ are less than or equal to $12$, and 
showed that the Turaev-Viro-Ocneanu invariant from the $E_6$-subfactor  
distinguishes the lens spaces $L(3,1)$ and $L(3,2)$. So, 
the Turaev-Viro-Ocneanu invariant distinguishes orientations for specific manifolds. 
That is remarkable result since the original Turaev-Viro invariant cannot distinguish orientations. 
 (More precisely, it coincides with the square of absolute value of the 
Reshetikhin-Turaev invariant \cite{Turaev}.) 
\par 
In \cite{KSW}, we showed that a Turaev-Viro-Ocneanu TQFT $Z$ from a subfactor has a Verlinde 
basis and that the invariant $Z(M)$ of a closed oriented $3$-manifold $M$ 
is given by the formula
$$Z(M)=\sum_{i_1,\dots,i_r=0}^m S_{0i_1} \dots S_{0i_r} J(L;i_1, \dots, i_r),$$
if $M$ is obtained from the $3$-sphere $S^3$ by Dehn surgery along a framed link $L = L_1 \cup \dots \cup L_r$.  
Here, $S_{ij}$ is the matrix component of $Z(S)$ 
with respect to a Verlinde basis $\{ v_i \}_{i=0}^m$ of $Z(S^1\times S^1)$, 
$S$ is the orientation preserving diffeomorphism on 
$S^1\times S^1$ corresponding to the modular transformation $S:\tau \longmapsto -1/\tau $, 
and $J(L;i_1,\cdots ,i_r)$ $(i_1,\cdots ,i_r=0,1,\cdots ,m)$ 
are framed link invariants of $L$ determined by the Turaev-Viro-Ocneanu invariant of the complement of $L$ in $S^3$. 
By the above formula, if we want to compute the Turaev-Viro-Ocneanu invariant 
of a closed 3-manifold $M$, we need to compute the $S$-matrix and the framed link 
invariants $J(L;i_1, \dots, i_r)$. 
Since we know that the method of a concrete construction of Verlinde basis (see \cite{KSW}), 
we can compute the $S$-matrix with respect to this Verlinde basis in principle. 
So, if a $3$-manifold $M$ is obtained from $S^3$ 
by Dehn surgery along an \lq\lq easy'' framed link $L$, then we can compute the 
$S$-matrix, the framed link invariants $J(L;i_1, \dots, i_r)$ and $Z(M)$. 
\par 
The concept of a tube algebra, which plays a crucial role in the Turaev-Viro-Ocneanu TQFT, 
was first introduced by Ocneanu \cite{Ocneanu1}. 
In analysis of the Longo-Rehren subfactor, which corresponds to 
the center construction in the sense of Drinfel'd, M. Izumi \cite{Izumi4}
formulated a tube algebra in terms of sectors for a finite closed 
system of endomorphisms of a type III subfactor. He also explicitly gave an action 
of $SL(2,\mathbb{Z})$ on the center of this tube algebra in the language 
of sectors, and derived several formulas on Turaev-Viro-Ocneanu invariants 
of lens spaces for concretely given subfactors \cite{Izumi5}. 
However, it is not clear how the 
$S$-matrix in Izumi's sector theory and the $S$-matrix in the Turaev-Viro-Ocneanu TQFT are related. 
\par 
In this paper, we establish a rigorous correspondence between the two tube algebras, 
that one comes from the Turaev-Viro-Ocneanu TQFT introduced by Ocneanu 
and another comes from the sector theory introduced by Izumi, 
and construct a canonical isomorphism between the centers of the two tube algebras, 
which is a conjugate linear isomorphism preserving the products of the two algebras 
and commuting with the actions of $SL(2,\mathbb{Z})$ (see Theorem \ref{Theorem1} and Corollary \ref{Corollary}). 
Moreover, we compute Turaev-Viro-Ocneanu invariants from several subfactors for basic $3$-manifolds including 
lens spaces and Brieskorn $3$-manifolds by using Izumi's data written in terms of sectors \cite{Izumi5} 
and using the Dehn surgery formula. 
One of the most important results on computations is that the $3$-sphere $S^3$ 
and the Poincar\'{e} homology $3$-sphere $\Sigma (2,3,5)$ are distinguished by the Turaev-Viro-Ocneanu invariant 
from the exotic subfactor constructed by Haagerup and Asaeda \cite{Haagerup, AsaedaHaagerup}, 
and $L(p,1)$ and $L(p,2)$ are distinguished by the Turaev-Viro-Ocneanu invariant   
from a generalized $E_6$-subfactor with $\mathbb{Z}/p\mathbb{Z}$ for $p=3,5$ \cite{Izumi5}. 
From this fact, it is natural for us to expect that the lens spaces $L(7,1)$ and $L(7,2)$ are distinguished 
by a generalized $E_6$-subfactor with $\mathbb{Z}/7\mathbb{Z}$. 
However, we do not know yet that there exists such a subfactor. 
\par 
This paper is organized as follows. 
In Section 2, we review the Turaev-Viro-Ocneanu topological quantum field theory from a subfactor. 
In Section 3, we discuss on the fusion algebra associated with the Turaev-Viro-Ocneanu TQFT \cite{KSW}. 
The product of the fusion algebra and the representation of $SL(2,\mathbb{Z})$ 
associated with the Turaev-Viro-Ocneanu TQFT are concretely calculated by using singular triangulations. 
In Section 4, we construct a conjugate linear isomorphism between Izumi's tube algebra and Ocneanu's one 
preserving the products of algebras and commuting with the actions of $SL(2,\mathbb{Z})$. 
In Section 5, we calculate Turaev-Viro-Ocneanu invariants from several subfactors for some basic $3$-manifolds 
including lens spaces and Brieskorn $3$-manifolds based on Izumi's data using the Dehn surgery formula, 
and derive formulas in terms of the $S$- and $T$-matrices. 
\par 
Throughout this paper, we assume that a closed oriented surface is geometrically realized as a polyhedron. 
We use the following notations. 
The $r$-simplex with vertices $v_0,v_1,\cdots ,v_r$ in an Euclidean space is denoted by $|v_0v_1\cdots v_r|$, and  
the simplicial complex consisting of all faces of $|v_0v_1\cdots v_r|$ is denoted by $K(|v_0v_1\cdots v_r|)$. 
For a $2$-simplex $|v_0v_1|$ we denote by $\langle v_0,v_1\rangle $ 
the oriented edge with the direction from $v_1$ to $v_0$. 
\par 
We refer to Evans and Kawahigashi's book \cite{EK} as a general reference on subfactor theory, 
and refer to Turaev's book \cite{Turaev} as a general reference on topological quantum field theory.

\par \bigskip 
\section{\kern-1em .\ \ Turaev-Viro-Ocneanu TQFT's from Subfactors} 
\par 
In  this section, we review the Turaev-Viro-Ocneanu TQFT's arising from subfactors \cite{EK, Ocneanu1} in terms of sectors. 
It is formulated by a similar method of Turaev and Viro \cite{TV} and by using initial data derived from subfactors. 
A theory of sectors in subfactor theory was first established by Longo \cite{Longo2} based on ideas 
in quantum field theory, and it was developed by Izumi \cite{Izumi1, Izumi3}. 
To describe the Turaev-Viro-Ocneanu TQFT's using sectors, we recall some notions in subfactor theory. 
\par 
Let $M$ be an infinite factor. We denote by $\text{End}(M)_0$ the set of $\ast $-endomorphisms 
$\rho :M\longrightarrow M$ such that the index of the subfactor $\rho (M)\subset M$ is finite. 
For $\rho ,\sigma \in \text{End}(M)_0$ the {\it intertwiner space} $(\rho ,\sigma )$ is a vector space over $\mathbb{C}$ 
defined by 
$$(\rho ,\sigma ):=\{ V\in M\ \vert \ V\rho (x)=\sigma (x)V\ \text{for}\ x\in M\} .$$
A $\ast $-endomorphism $\rho \in \text{End}(M)_0$ is called {\it irreducible}, if $(\rho ,\rho )=\mathbb{C}id_M$. 
If $\rho $ is irreducible, then for any $\sigma \in \text{End}(M)_0$ 
the intertwiner space $(\rho ,\sigma )$ is a Hilbert space with the inner product defined by 
$$\langle V, W\rangle =W^{\ast }V,\qquad V, W\in (\rho ,\sigma ).$$
For $\rho \in \text{End}(M)_0$ we denote by $d(\rho )$ the square root of the minimal index of $M\supset \rho (M)$, 
and call it the statistical dimension of $\rho $. 
It is known that for every $\rho \in \text{End}(M)_0$ 
there exists a $\ast $-endomorphism $\bar{\rho} \in \text{End}(M)_0$ 
and a pair of intertwiners $R_{\rho }\in (id, \bar{\rho }\rho )$, 
$\overline{R}_{\rho }\in (id, \rho \bar{\rho })$ such that 
$$\overline{R}_{\rho }^{\ast}\rho (R_{\rho })=R_{\rho }^{\ast}\bar{\rho }(\overline{R}_{\rho })
=\dfrac{1}{d(\rho )},\ R_{\rho }^{\ast }R_{\rho }=\overline{R}_{\rho }^{\ast }\overline{R}_{\rho }=1.$$

Let $\rho _1$ and $\rho _2$ be two elements in $\text{End}(M)_0$. 
We say that $\rho _1$ and $\rho _2$ are equivalent if there exists a unitary $u\in M$ such that $u\rho _1(x)=\rho _2(x)u$ for all $x\in M$. 
We denote by $\text{Sect}(M)$ the equivalence classes of $\text{End}(M)_0$, and denote by $[\rho ]$ the class of 
$\rho \in \text{End}(M)_0$ in $\text{Sect}(M)$. 
Each element in $\text{Sect}(M)$ is called a {\it sector} of $M$. 
The set $\text{Sect}(M)$ becomes a $\ast $-semiring over $\mathbb{C}$ 
with the sum $[\rho ]\oplus [\sigma ]$, 
the product $[\rho ][\sigma ]=[\rho \circ \sigma ]$ 
and the conjugation $\overline{[\rho ]}=[\bar{\rho}]$ \cite{Izumi1, Longo2}. 

\par 
A finite subset $\Delta =\{ \rho _0,\rho _1,\cdots ,\rho _n\}$ of $\text{End}(M)_0$ is called 
{\it a finite system of $\text{End}(M)_0$} closed under sector operations 
if the following four conditions are satisfied \cite{Izumi4}.

\begin{enumerate}
\item[(i)] $[\rho _i]=[\rho _j]$ if and only if $i=j$.
\item[(ii)] $\rho_0=id_M$.
\item[(iii)] for all $i$ there exists $j$ such that $\overline{[\rho _i]}=[\rho _j]$. 
\item[(iv)] there exist non-negative integers $N_{ij}^k$  such that 
$$[\rho _i][\rho _j]=\textstyle\bigoplus\limits_{k=0}^nN_{ij}^k[\rho _k].$$
\end{enumerate}

We note that the condition (iv) is equivalent to the following condition. 

\begin{enumerate}
\item[(iv)$'$] for all $i,j,k=0,1,\cdots ,n$, there exists an orthonormal basis  
$\{ (T_{ij}^k)_{\nu }\} _{\nu =0}^{N_{ij}^k}$ of $(\rho _k,\rho _i\rho _j)$ such that 
$$\sum\limits_{k=0}^n\sum\limits_{\nu =1}^{N_{ij}^k}(T_{ij}^k)_{\nu }(T_{ij}^k)_{\nu }^{\ast }=1,\ 
\rho _i\rho _j=\sum\limits_{k=0}^n\sum\limits_{l=1}^{N_{ij}^k}(T_{ij}^k)_{\nu }\rho _k(T_{ij}^k)_{\nu }^{\ast }.$$
\end{enumerate}

If $N\subset M$ is an inclusion of infinite factors with finite index and finite depth, 
then we can obtain a finite system of $\text{End}(M)_0$ closed under sector operations 
by taking irreducible components in $([\iota ] \overline{[\iota ]})^n$, $([\iota ] \overline{[\iota ]})^n 
[\iota ] $, $\overline{[\iota ]} ([\iota ] \overline{[\iota ]})^n $ and 
$(\overline{[\iota ]} [\iota ])^n$, where $\iota : N \hookrightarrow M$ is the inclusion 
map (See \cite{KSW} for detail). 
We call it a {\it finite system of $\text{End}(M)_0$ obtained from $N\subset M$}. 
\par 
Frobenius reciprocities for sectors were established by Izumi \cite{Izumi3} as an analogue to Frobenius reciprocities  
for group representations and bimodules with left and right actions of II$_1$-factors \cite{EK}. 
Let $\Delta $ be a finite system of $\text{End}(M)_0$ obtained from a subfactor $N\subset M$ 
of an infinite factor $M$ with finite index and finite depth. 
For $\rho ,\eta ,\zeta \in \Delta $, let $\mathcal{H}_{\rho \eta }^{\zeta }$ 
denote the intertwiner space $(\zeta , \rho \eta )$. 
Then, two conjugate linear maps 
$(\widetilde{\mathstrut  \cdot })_{\rho \eta }^{\zeta }:\mathcal{H}_{\rho \eta }^{\zeta }\longrightarrow \mathcal{H}_{\bar{\rho } \zeta }^{\eta }$ 
and $(\widehat{\mathstrut \cdot })_{\rho \eta }^{\zeta }:\mathcal{H}_{\rho \eta }^{\zeta }\longrightarrow \mathcal{H}_{\zeta \bar{\eta} }^{\rho }$ 
are defined by 
$$(\widetilde{A})_{\rho \eta }^{\zeta }=\sqrt{\dfrac{d(\rho )d(\eta )}{d(\zeta )}}\bar{\rho }(A^{\ast })R_{\rho },\quad 
(\widehat{A})_{\rho \eta }^{\zeta }=\sqrt{\dfrac{d(\rho )d(\eta )}{d(\zeta )}}A^{\ast} \rho (\bar{R}_{\eta })$$
for all $A\in \mathcal{H}_{\rho \eta }^{\zeta }$, respectively. 
They give rise to an action of the symmetric group $S_3$ of degree $3$ between the following spaces:
$$\mathcal{H}_{\rho \eta }^{\zeta }, \mathcal{H}_{\bar{\zeta} \rho }^{\bar{\eta }}, 
\mathcal{H}_{\eta \bar{\zeta }}^{\bar{\rho } }, \mathcal{H}_{\bar{\eta} \bar{\rho} }^{\bar{\zeta }}, 
\mathcal{H}_{\bar{\rho } \zeta }^{\eta }, \mathcal{H}_{\zeta \bar{\eta }}^{\rho }.$$ 
The maps $(\widetilde{\mathstrut \cdot })_{\rho \eta }^{\zeta }$ 
and $(\widehat{\mathstrut \cdot })_{\rho \eta }^{\zeta }$ are called the (left) Frobenius reciprocity maps \cite{Izumi3}. 
For simplicity, we frequently write $\widetilde{A}$ and $\widehat{A}$ instead of 
$(\widetilde{A})_{\rho \eta }^{\zeta }$ and $(\widehat{A})_{\rho \eta }^{\zeta }$, respectively. 

\par \bigskip 
We will describe the Turaev-Viro-Ocneanu TQFT's arising from subfactors in the setting of sectors \cite{EK, Ocneanu1}. 
Let $\Delta $ be a finite system of $\text{End}(M)_0$ obtained from a subfactor $N\subset M$ 
of an infinite factor $M$ with finite index and finite depth. 
We fix an orthonormal basis $\mathcal{B}_{\rho  \eta }^{\zeta }$ of $\mathcal{H}_{\rho  \eta }^{\zeta }$ 
for each $\rho , \eta , \zeta \in \Delta $.  
\par 
A simplicial complex $\mathcal{K}$ which edges are oriented is called {\it a locally ordered complex}, 
if there is no cyclic order for the $3$ edges of any triangle in $\mathcal{K}$. 
By a {\it color} of a locally ordered complex $\mathcal{K}$, we mean a map 
$$\xi : \{ \text{the oriented edges and the triangles in $\mathcal{K}$ }\} 
\longrightarrow \Delta \cup \bigcup\limits_{\rho ,\eta ,\zeta \in \Delta } \mathcal{B}_{\rho  \eta }^{\zeta }$$
satisfying the following two conditions.
\begin{enumerate}
\item[(i)] $\xi (E)\in \Delta $ for each oriented edge $E$ in $\mathcal{K}$.
\item[(ii)] If $\xi (\langle v_0,v_1\rangle )=\rho ,\ \xi (\langle v_1,v_2\rangle )=\eta ,\ 
\xi (\langle v_0,v_2\rangle )=\zeta $ for a triangle $|v_0v_1v_2|$ in $\mathcal{K}$, 
then $\xi (\vert v_0v_1v_2\vert )\in \mathcal{B}_{\rho  \eta }^{\zeta }$. 
\end{enumerate}

We follow the convention that 
$\xi (-E)=\overline{\xi (E)}$ for each oriented edge $E$ in $\mathcal{K}$, 
where $-E$ denotes the same edge $E$ with opposite orientation.

For a color $\xi $ we define a positive real number $d(\xi )$ by 
$$d(\xi )=\prod _{e\text{ : edges of $K$}}d(\xi (e)).$$ 

\begin{figure}[hbtp]
\setlength{\unitlength}{1cm}
\begin{center}
\vspace{-0.3cm}
\includegraphics[height=3cm]{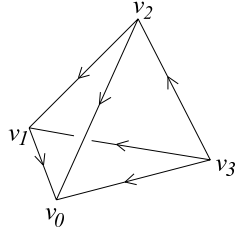}
\caption{a tetrahedron \label{Figure1}}
\end{center}
\end{figure}

\begin{figure}[hbtp]
\setlength{\unitlength}{1cm}
\begin{center}
\vspace{-0.3cm}
\includegraphics[height=3cm]{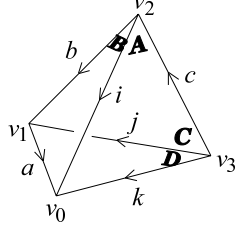}
\caption{a colored tetrahedron \label{Figure1a}}
\end{center}
\end{figure}

Let $X$ be a compact oriented $3$-manifold possibly with boundary. 
We fix a locally ordered complex $\mathcal{K}$ that is a triangulation of $\partial X$, 
and choose a locally ordered complex $\mathcal{T}$ that is a triangulation of $X$ 
satisfying $\partial \mathcal{T}=\mathcal{K}$. 
To a color $\varphi $ of $\mathcal{T}$ and a tetrahedron $\sigma =\vert
v_0v_1v_2v_3\vert $ in $\mathcal{T}$ depicted as in Figure \ref{Figure1}, 
we assign a complex number defined by 

$$\dfrac{1}{\sqrt{d(i)d(j)}}A^{\ast }B^{\ast }a(C)D, $$
where $\varphi (\langle v_0,v_1\rangle )=a$, $\varphi (\langle v_1,v_2\rangle )=b$, 
$\varphi (\langle v_2,v_3\rangle )=c$, $\varphi (\langle v_0,v_2\rangle )=i$, 
$\varphi (\langle v_1,v_3\rangle )=j$, $\varphi (\langle v_0,v_3\rangle )=k$, 
$\varphi (|v_1v_2v_3|)=A$, $\varphi (|v_0v_1v_2|)=B$, $\varphi (|v_1v_2v_3|)=C$, 
$\varphi (|v_0v_1v_3|)=D$ (see Figure \ref{Figure1a}). 
We denote the above complex number or its complex conjugate by $W(\sigma ; \varphi  )$ according to compatibility of 
orientations for $M$ and $\sigma$. 
Here, the orientation for $\sigma $ is given by the order $v_0<v_1<v_2<v_3$. 

\par 
For a color $\xi $ of $\partial \mathcal{T}$, we set 

$$Z(X;\mathcal{T},\xi ):=\lambda ^{-\sharp \mathcal{T}^{(0)}+\frac{\sharp (\partial \mathcal{T})^{(0)}}{2}}\sqrt{d(\xi )}
\sum \limits_{\substack{\varphi \text{:color of $\mathcal{T}$}\\ \varphi \vert _{\partial \mathcal{T}}=\xi }}\kern-0.7em
d(\varphi \vert _{\mathcal{T}-\partial \mathcal{T}})\kern-0.5em \prod \limits_{\sigma :\text{tetrahedra}}W(\sigma ;\varphi ),$$
where $\mathcal{T}^{(0)}$ and $\partial \mathcal{T}^{(0)}$ are the sets of vertices of $\mathcal{T}$ 
and $\partial \mathcal{T}$, respectively, and $\lambda =\sum\limits_{i=0}^nd(\rho _i)^2$, 
which is called the {\it global index} of $\Delta $. 
It is proved that 
$Z(X;\mathcal{T},\xi )$ does not depend on the choice of orthonormal bases $\mathcal{B}_{\rho  \eta }^{\zeta }$\ 
$(\rho ,\eta ,\zeta \in \Delta )$, 
and the choice of triangulations $\mathcal{T}$ of $X$ such that $\partial \mathcal{T}=\mathcal{K}$ \cite{EK}. 
In particular, if $X$ is closed, 
then $Z(X;\mathcal{T},\xi )$ is a topological invariant of the closed oriented $3$-manifold $X$. So we may denote it by $Z(X)$, 
and refer this invariant as the Turaev-Viro-Ocneanu invariant of $X$. 

\par \bigskip 
The method of the construction of $Z(X;\mathcal{T},\xi )$ as mentioned above gives rise to a functor
$Z^{\Delta }$ that satisfies the axioms of topological quantum field theory as posed by
Atiyah \cite{Atiyah}. (Although the same functor was denoted by $Z_{\Delta }$ in the previous paper \cite{KSW}, 
we use the notation $Z^{\Delta }$ for such a functor in order to improve the appearance of notations used later.) 
We briefly describe the method of constructing such a functor $Z^{\Delta }$ following Turaev and Viro's paper \cite{TV} 
and Yetter's paper \cite{Yetter}.
\par 
Let $\mathcal{P}(\Sigma )$ denote the set of locally ordered complexes 
that are triangulations of a closed oriented surface $\Sigma $.  
For $\mathcal{K}\in \mathcal{P}(\Sigma )$, we denote by $V(\Sigma ;\mathcal{K})$ 
the $\mathbb{C}$-vector space freely spanned by the colors of $\mathcal{K}$. If $\Sigma =\emptyset $, 
then we set $V(\Sigma ; \mathcal{K}):=\mathbb{C}$. 
By a $(2+1)$-dimensional cobordism with triangulated boundary 
we mean a $(2+1)$-dimensional cobordism $(X;\Sigma _1,\Sigma _2)$ 
between oriented closed surfaces $\Sigma _1$ and $\Sigma _2$ triangulated by locally ordered complexes 
$\mathcal{K}_1$ and $\mathcal{K}_2$, respectively. 
We denote it by $(X;\mathcal{K}_1,\mathcal{K}_2)$. 
For such a cobordism $W=(X;\mathcal{K}_1,\mathcal{K}_2)$, 
a $\mathbb{C}$-linear map $\Phi _W:V(\Sigma _1;\mathcal{K}_1)\longrightarrow V(\Sigma _2;\mathcal{K}_2)$ is defined by 
$$\Phi _W(\xi _1)=\sum\limits_{\xi _2\text{ : the colors of $\mathcal{K}_2$}} Z(X;\mathcal{T}, \xi _1\cup \xi _2)\xi _2,$$
for all colors $\xi _1$ of $\mathcal{K}_1$, where 
$\mathcal{T}$ is a triangulation of $X$ such that $\partial \mathcal{T}=\mathcal{K}_1\cup \mathcal{K}_2$, and 
$\xi _1\cup \xi _2$ is the color of $\mathcal{T}$ given by $(\xi _1\cup \xi _2)\vert _{\mathcal{K}_i}=\xi _i\ (i=1,2)$. 
\par 
For a locally ordered complex $\mathcal{K}\in \mathcal{P}(\Sigma )$ 
and a stellar subdivision $\mathcal{K}'$ of $\mathcal{K}$, which is denoted by $\mathcal{K}'\leq \mathcal{K}$, 
we have a cobordism with triangulated boundary 
$Id_{\mathcal{K}',\mathcal{K}}:=(\Sigma \times [0,1];\mathcal{K}',\mathcal{K})$ such that $\Sigma \times \{ 0\} $ and 
$\Sigma \times \{ 1\} $ are triangulated by $\mathcal{K}'$ and $\mathcal{K}$, respectively. 
This cobordism induces the $\mathbb{C}$-linear map 
$\Phi _{Id_{\mathcal{K}',\mathcal{K}}}:V(\Sigma ;\mathcal{K})\longrightarrow V(\Sigma ;\mathcal{K}')$. 
Let $\iota _{\mathcal{K}}:V(\Sigma ;\mathcal{K})\longrightarrow \bigoplus_{\mathcal{K}\in \mathcal{P}(\Sigma )}V(\Sigma ;\mathcal{K})$ be the canonical injection, and 
$W(\Sigma )$ the subspace of $\bigoplus_{\mathcal{K}\in \mathcal{P}(\Sigma )}V(\Sigma ;\mathcal{K})$ spanned by 
$$\{ \iota _{\mathcal{K}'}(x)-(\iota _{\mathcal{K}}\circ \Phi _{\mathcal{K}',\mathcal{K}})(x)\ \vert \ x\in V(\Sigma ;\mathcal{K}'),\ \mathcal{K},\mathcal{K}'\in \mathcal{P}(\Sigma ),\ \mathcal{K'}\leq \mathcal{K}\} .$$
We consider the quotient space 
$$Z^{\Delta }(\Sigma )=\bigoplus _{\mathcal{K}\in \mathcal{P}(\Sigma )}V(\Sigma ;\mathcal{K})/W(\Sigma ) .$$ 
Since $\Phi _{Id_{\mathcal{K},\mathcal{K}}}^2=\Phi _{Id_{\mathcal{K},\mathcal{K}}}$, the following lemma holds. 

\par \bigskip 
\begin{lem}[Yetter \cite{Yetter}]
\label{Lemma1}  
Let $u_{\mathcal{K}}:V(\Sigma ;\mathcal{K})\longrightarrow Z^{\Delta }(\Sigma )$ be the composition of 
$\iota _{\mathcal{K}}$ and the natural projection 
$\bigoplus _{\mathcal{K}\in \mathcal{P}(\Sigma )}V(\Sigma ;\mathcal{K})\longrightarrow Z^{\Delta }(\Sigma )$. 
Then, $u_{\mathcal{K}}$ is surjective and 
$\text{Ker}\kern0.2em u_{\mathcal{K}}=\text{Ker}\kern0.2em \Phi _{Id_{\mathcal{K},\mathcal{K}}}$. 
In particular, $Z^{\Delta }(\Sigma )$ is finite-dimensional, and 
$Z^{\Delta }(\Sigma )\cong V(\Sigma ;\mathcal{K})/\text{Ker}\kern0.2em \Phi _{Id_{\mathcal{K},\mathcal{K}}}$ 
as vector spaces. 
\end{lem}

By the above lemma, we can verify that 
the linear map $\Phi _W$ induces a linear map 
$Z_W^{\Delta }: Z^{\Delta }(\Sigma _1)\longrightarrow Z^{\Delta }(\Sigma _2)$ 
for each cobordism with triangulated boundary $W$, and $Z_W^{\Delta }$ does not depend on the choice of triangulations. 
\par 
Let $f:(\Sigma _1, \mathcal{K}_1) \longrightarrow (\Sigma _2, \mathcal{K}_2)$ be an isomorphism 
between closed oriented surfaces $\Sigma _1$  and $\Sigma _2$ 
triangulated by locally ordered complexes $\mathcal{K}_1$ and $\mathcal{K}_2$, respectively. 
Then, a $\mathbb{C}$-linear isomorphism 
$\Phi (f):V(\Sigma _1;\mathcal{K}_1)\longrightarrow V(\Sigma _2;\mathcal{K}_2)$ 
is defined by
$$\Phi (f)(\xi _1)=\xi _1\circ f^{-1},$$
for all colors $\xi _1$ of $\mathcal{K}_1$. 
Here, $\xi _1\circ f^{-1}$ is a color of $\mathcal{K}_2$ such that 
$(\xi _1\circ f^{-1})(E)=\xi _1(f^{-1}(E))$ for each edge and face $E$ of $\mathcal{K}_2$. 
\par 
Let $f:\Sigma _1 \longrightarrow \Sigma _2$ be an 
orientation preserving PL-homeomorphism between closed oriented surfaces. 
For $\mathcal{K}_i\in \mathcal{P}(\Sigma _i)\ (i=1,2)$, 
there exists a refinement $\mathcal{K}_i'$ of $\mathcal{K}_i$ such that $f$ is 
an isomorphism from $\mathcal{K}_1'$ to $\mathcal{K}_2'$ as a locally ordered complex.
Then, we have a linear map
$$f_{\#}:=\Phi _{Id_{\mathcal{K}_2',\mathcal{K}_2}}\circ \Phi (f)\circ \Phi_{Id_{\mathcal{K}_1,\mathcal{K}_1'}}.$$

By Lemma \ref{Lemma1}, we can verify that 
the linear map $f_{\#}$ induces the linear isomorphism 
$Z^{\Delta }(f): Z^{\Delta }(\Sigma _1)\longrightarrow Z^{\Delta }(\Sigma _2)$, and 
$Z^{\Delta }(f)$ does not depend on the choice of triangulations. 
\par 
In the above manner, 
we have a functor $Z^{\Delta }$ that assigns 
to each closed oriented surface $\Sigma $ the $\mathbb{C}$-vector space $Z^{\Delta }(\Sigma )$, 
to each $(2+1)$-dimensional cobordism $W$ the $\mathbb{C}$-linear map $Z_W^{\Delta }$, and 
to each orientation preserving PL-homeomorphism $f$ 
between closed oriented surfaces the $\mathbb{C}$-linear isomorphism $Z^{\Delta }(f)$. 
This functor $Z^{\Delta }$ satisfies the axioms for $(2+1)$-dimensional TQFT. 
Furthermore, the TQFT $Z^{\Delta }$ is unitary, since 
$Z^{\Delta }(\Sigma )$ becomes a Hilbert space with the inner product induced from the hermitian form on $V(\Sigma ;\mathcal{K})$, which is defined by  
$$\langle \xi _0, \xi _1 \rangle _{TQFT}=Z^{\Delta }(\Sigma \times [0,1];\mathcal{K}\times [0,1], \xi _0\cup \xi _1)$$ 
for all pairs of colors $\xi _0$ and $\xi _1$ of $\mathcal{K}$ \cite{SuzukiWakui, Turaev}.

\par \bigskip 
\section{\kern-1em .\ \ Fusion Algebras Associated with Turaev-Viro-Ocneanu TQFT's}
\par 
For a short while, we consider a $(2+1)$-dimensional TQFT $Z$ in a general situation apart from subfactors. 
We describe the definition of the fusion algebra in the $(2+1)$-dimensional TQFT $Z$, and 
introduce our version of Verlinde basis of $Z(S^1\times S^1)$ (see Definition \ref{Definition}). 
\par
\begin{figure}[hbtp]
\begin{center}
\includegraphics[height=3cm]{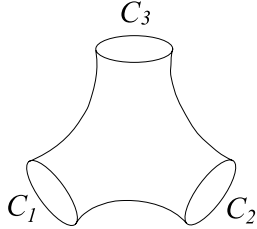}
\caption{the compact oriented surface $Y$ \label{Figure2}}
\end{center}
\end{figure}

 
Let $W$ be the cobordism $(Y\times S^1; \Sigma _1\sqcup \Sigma _2, \Sigma _3)$, where $Y$ is
the $3$-holed sphere in $\mathbb{R}^3$
depicted in Figure \ref{Figure2} and $\Sigma _i=C_i\times S^1$ for $i=1,2,3$. 
Then, $W$ induces a 
linear map $Z_W: Z(S^1\times S^1)\otimes Z(S^1\times S^1)\longrightarrow
Z(S^1\times S^1)$. It 
can be easily verified that the map $Z_W$ gives an associative algebra structure
on $Z(S^1\times S^1)$. The identity element of this algebra is given by
$Z_{W_0}(1)$, where $W_0:=(D^2\times S^1; \emptyset , S^1\times S^1)$. We
call this algebra the {\it fusion algebra associated with $Z$}.

\begin{figure}[hbtp]
\begin{center}
\setlength{\unitlength}{1cm}
\includegraphics[height=6cm]{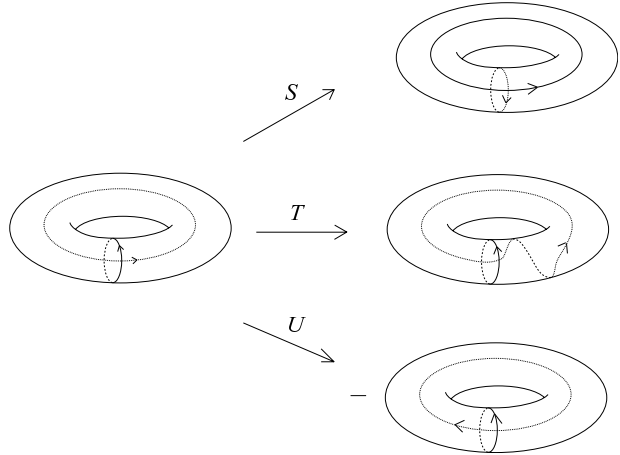}
\caption{the action of $SL(2,\mathbb{Z})$ \label{Figure3}}
\end{center}
\end{figure}

Let us recall that the mapping class group $\varGamma_{S^1\times S^1}$ of the torus $S^1\times S^1$ 
is isomorphic to the group 
$SL(2,\mathbb{Z})$ of integral $2\times 2$-matrices with determinant $1$. 
It is well-known that this group is generated by $S=\begin{pmatrix} 0 & 1 \\ -1 & 0
\end{pmatrix}$ and $T=\begin{pmatrix} 1 & 0 \\ 1 & 1\end{pmatrix}$ with relations
$S^4=I,\ (ST)^3=S^2$. The matrices $S$ and $T$ correspond to the orientation preserving diffeomorphisms 
from $S^1\times S^1$ to $S^1\times S^1$ which are defined by 

\begin{equation}
S(z,w)=(\bar{w},z), \qquad 
T(z,w)=(zw,w) \label{eq1}
\end{equation}

\noindent 
for all $(z,w)\in S^1\times S^1$, respectively, where 
we regard $S^1$ as the set of complex numbers of absolute value $1$. 
To define the Verlinde basis, 
we need one more orientation preserving diffeomorphism 
$U: S^1\times S^1\longrightarrow - S^1\times S^1$ 
defined by $U(z,w)=(z,\bar{w})$ for all $(z,w)\in S^1\times S^1$ (see Figure \ref{Figure3}). 

\par \bigskip 
\begin{defn}[\cite{KSW}]\ \ 
\label{Definition}
Let $Z$ be a $(2+1)$-dimensional TQFT. A basis $\{ v_i\} _{i=0}^m$ of $Z(S^1\times S^1)$ 
is said to be a {\it Verlinde basis} if it has the following properties. 
\begin{enumerate} 
\item[(i)] $v_0$ is the identity element of the fusion algebra associated with $Z$. 
\item[(ii)] 
\begin{enumerate} 
\item[(a)] $Z(S)$ is represented by a unitary and symmetric matrix with respect to the basis $\{ v_i\} _{i=0}^m$.
\item[(b)] $Z(S)^2v_0=v_0$, and $Z(S)^2v_i\in \{ v_j\} _{j=0}^m$ for all $i$. 
\item[(c)] We define $S_{ji}\in \mathbb{C}$ by $Z(S)v_i=\sum_{i=0}^mS_{ji}v_j$. Then, 
\begin{enumerate} 
\item[1.] $S_{i0}\not= 0$ for all $i$. 
\item[2.] $N_{ij}^k:=\sum_{l=0}^m\frac{S_{il}S_{jl}\overline{S_{lk}}}{S_{0l}}\ (i,j,k=0,1,\cdots ,m)$ 
coincide with the structure constants of the fusion algebra with respect to $\{ v_i\} _{i=0}^m$. 
\end{enumerate}
\end{enumerate}
\item[(iii)] $Z(T)$ is represented by a diagonal matrix with respect to the basis $\{ v_i\} _{i=0}^m$. 
\item[(iv)] $Z(U)v_i=v_{i}^{\ast }$ for all $i$ 
under the identification $Z(-S^1\times S^1)\cong Z(S^1\times S^1)^{\ast }$. 
\end{enumerate}
\end{defn}

\par \bigskip 
\begin{rems}
\label{remarks}
1. If $\{ v_i\} _{i=0}^m$ is a Verlinde basis, then $S_{i0}$ is a real number for all $i$.  
\par 
\noindent 
2. A Verlinde basis is unique up to order of elements, 
since $w_i=S_{0i}\sum_{j=0}^m\overline{S_{ji}}v_j$\ $(i=0,1,\cdots ,m)$ 
are all orthogonal primitive idempotents in the fusion algebra satisfying $1=w_0+w_1+\cdots +w_m$. 
This fact follows from that the $S$-matrix diagonalizes the fusion rules in conformal field theory \cite{Verlinde}. 
\par 
\noindent 
3. The map $\overline{\mathstrut \ \cdot \ }:\{ 0,1,\cdots ,m\} \longrightarrow \{ 0,1,\cdots ,m\} $ defined by $Z(S)^2v_i=v_{\bar{i}}$ is an involution satisfying $\bar{0}=0$. 
\par 
\noindent 
4. The last condition (iv) was introduced in \cite{Wakui} and modified in \cite{SuzukiWakui}. 
\end{rems}

\par \bigskip 
We go back to the setting of subfactors. 
Let $\Delta $ be a finite system of $\text{End}(M)_0$ obtained from a subfactor $N\subset M$ 
of an infinite factor $M$ with finite index and finite depth. 
To calculate the product of the fusion algebra and the representation of $SL(2,\mathbb{Z})$ 
associated with the Turaev-Viro-Ocneanu TQFT $Z^{\Delta }$ arising from subfactors, 
we show that the fusion algebra associated with $Z^{\Delta }$ 
can be regarded as a subalgebra of the tube algebra, 
which plays a crucial role in the Turaev-Viro-Ocneanu TQFT. 
We review the definition of a tube algebra \cite{KSW} in $Z^{\Delta }$, 
which was first introduced by Ocneanu \cite{Ocneanu1}. 
\par
\begin{figure}[hbt]
\setlength{\unitlength}{1cm}
\begin{center}
\includegraphics[height=4.5cm]{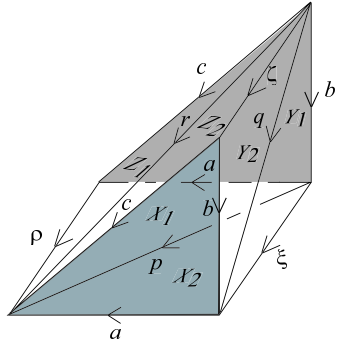}
\caption{the coefficient of the product of $\text{Tube}\Delta $ \label{Figure4}}
\end{center}
\end{figure}

 
The tube algebra $\text{Tube}\Delta $ is a finite-dimensional $C^{\ast }$-algebra over $\mathbb{C}$  
which is defined by 
$$\text{Tube}\Delta 
=\bigoplus\limits_{\rho ,\xi ,\zeta ,\eta \in \Delta }
\mathcal{H}_{\rho \eta }^{\zeta }\otimes \mathcal{H}_{\eta \xi }^{\zeta }$$
as $\mathbb{C}$-vector spaces. 
The product of $\text{Tube}\Delta $ is given by 
$$(X_1\otimes X_2)\cdot (Y_1\otimes Y_2)
=\delta _{\xi ,\eta }\dfrac{\lambda }{\sqrt{d(\rho )d(\zeta )}d(\xi )} 
\sum\limits_{\substack{r,c\in \Delta \\ Z_1\in \mathcal{B}_{\rho c}^r\\ Z_2\in \mathcal{B}_{c \zeta }^r}}
Z(D^2\times S^1;{}_{\rho ,a,p,X_1,X_2;\eta ,b,q,Y_1,Y_2}^{\zeta ,c,r,Z_1,Z_2})Z_1\otimes Z_2$$
for $X_1\otimes X_2\in \mathcal{H}_{\rho  a}^p\otimes \mathcal{H}_{a \xi }^p$, 
$Y_1\otimes Y_2\in \mathcal{H}_{\eta  b}^q\otimes \mathcal{H}_{b \zeta }^q$, 
where $\delta _{\xi ,\eta }$ is Kronecker's delta, 
$\lambda $ is the global index of $\Delta $, $\mathcal{B}_{\rho c}^r$ and 
$\mathcal{B}_{c \zeta }^r$ are orthonormal bases of $\mathcal{H}_{\rho c}^r$ and 
$\mathcal{H}_{c \zeta }^r$, respectively, 
and $Z(D^2\times S^1;{}_{\rho ,a,p,X_1,X_2;\eta ,b,q,Y_1,Y_2}^{\zeta ,c,r,Z_1,Z_2})$ is 
the Turaev-Viro-Ocneanu invariant of the solid torus $D^2\times S^1$ 
triangulated by $\mathcal{T}$ with colored boundary illustrated as in Figure \ref{Figure4}. 
(Here, the two triangles shaded inside are identified.) 
We remark that the normalization of the product is slightly different from one in Ocneanu's definition \cite{Ocneanu1}. 

\begin{figure}[htbp]
\setlength{\unitlength}{0.7cm}
\begin{center}
\begin{picture}(5,4)
\put(2,1){\framebox(2.5,2.5)}
\put(2,1){\line(1,1){2.5}}
\put(3.2,2.2){\line(1,0){0.2}}
\put(3.2,2.2){\line(0,1){0.2}}
\put(0.5,2){$\mathcal{K}$ =}
\put(3,0.85){$\ll $}
\put(3,3.38){$\ll $}
\put(1.8,2.2){$\vee $}
\put(4.32,2.2){$\vee $}
\end{picture}
\caption{a distinguished triangulation $\mathcal{K}$ of $S^1\times S^1$ \label{Figure5}}
\end{center}
\end{figure}
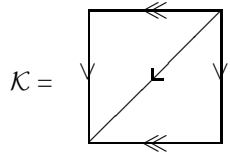

Let $V^{\Delta }(S^1\times S^1)$ denote the subalgebra of $\text{Tube}\Delta $ defined by 
$$V^{\Delta }(S^1\times S^1)=\bigoplus\limits_{\rho ,\zeta ,\eta \in \Delta }\mathcal{H}_{\rho \eta }^{\zeta }\otimes \mathcal{H}_{\eta \rho }^{\zeta }.$$
We regard $S^1\times S^1$ as a quotient space obtained from a square by identifying the opposite sides. 
Then, $V^{\Delta }(S^1\times S^1)$ is canonically identified with $V(S^1\times S^1;\mathcal{K})$ 
as a vector space, where $\mathcal{K}$ is the locally ordered complex depicted in Figure \ref{Figure5} 
that gives a singular triangulation of $S^1\times S^1$. 

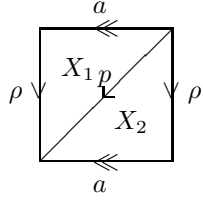
\begin{figure}[htbp]
\setlength{\unitlength}{0.7cm}
\begin{center}
\begin{picture}(5,4)
\put(2,1){\framebox(2.5,2.5)}
\put(2,1){\line(1,1){2.5}}
\put(3.2,2.2){\line(1,0){0.2}}
\put(3.2,2.2){\line(0,1){0.2}}
\put(3,0.4){$a$}
\put(3,3.8){$a$}
\put(3,0.85){$\ll $}
\put(3,3.38){$\ll $}
\put(1.4,2.2){$\rho $}
\put(4.8,2.2){$\rho $}
\put(1.8,2.2){$\vee $}
\put(4.32,2.2){$\vee $}
\put(3.1,2.5){$p$}
\put(2.4,2.6){$X_1$}
\put(3.4,1.6){$X_2$}
\end{picture}
\caption{a color of $\mathcal{K}$ \label{Figure6}}
\end{center}
\end{figure}

Under the identification $V^{\Delta }(S^1\times S^1)=V(S^1\times S^1;\mathcal{K})$, 
each element $X_1\otimes X_2$ in $V^{\Delta }(S^1\times S^1)$ such that  
$X_1\in \mathcal{H}_{\rho a}^p$ and $X_2\in \mathcal{H}_{a\rho }^p$, 
corresponds to the color of $\mathcal{K}$ depicted in Figure \ref{Figure6}. 

Now, we have two products on $V^{\Delta }(S^1\times S^1)$. 
One is obtained from the product of the tube algebra $\text{Tube}\Delta $ by restricting to $V^{\Delta }(S^1\times S^1)$. Another comes from the product of the fusion algebra $Z^{\Delta }(S^1\times S^1)$. They are not same, 
but are closely related as follows.

\par \bigskip 
\begin{prop}
\label{Proposition1}  
Let $\Delta $ be a finite system of $\text{End}(M)_0$ obtained from a subfactor $N\subset M$ of an infinite factor $M$ with finite index and finite depth. 
Let $P:V^{\Delta }(S^1\times S^1)\longrightarrow V^{\Delta }(S^1\times S^1)$ be the conjugate linear map defined by $P(X_1\otimes X_2)=X_2\otimes X_1$ for all $X_1\in \mathcal{H}_{\rho a}^p$ and $X_2\in \mathcal{H}_{a\rho }^p$. 
Then, the product $m_{\text{tube}}$ of the subalgebra $V^{\Delta }(S^1\times S^1)$ of $\text{Tube}\Delta $ 
and the product $m_{\text{fusion}}$ of the fusion algebra associated to $Z^{\Delta }$ 
are related by the following commutative diagram. 

\begin{equation*}
\begin{CD}
V^{\Delta }(S^1\times S^1) \otimes V^{\Delta }(S^1\times S^1)  @>m_{\text{tube}}>> V^{\Delta }(S^1\times S^1) \\ 
@VP\otimes PVV @VVPV \\ 
V^{\Delta }(S^1\times S^1) \otimes V^{\Delta }(S^1\times S^1)  @. V^{\Delta }(S^1\times S^1) \\ 
@Vu\otimes uVV @VVuV \\ 
Z^{\Delta }(S^1\times S^1)\otimes Z^{\Delta }(S^1\times S^1) @>m_{\text{fusion}}>> Z^{\Delta }(S^1\times S^1) 
\end{CD}
\end{equation*}
Here, $u$ is the universal arrow associated to $Z^{\Delta }(S^1\times S^1)$ defined as in Lemma \ref{Lemma1}. 
\end{prop}

\begin{figure}[hbtp]
\setlength{\unitlength}{1cm}
\begin{center}
\includegraphics[height=4cm]{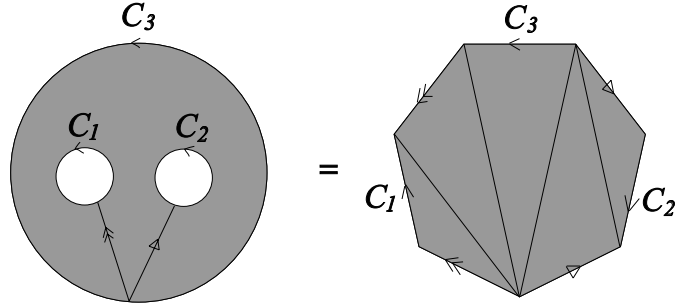}
\caption{a triangulation of $Y$ \label{Figure7}}
\end{center}
\end{figure}

\indent 
{\bf Proof.} 
The $3$-holed sphere $Y$ in $\mathbb{R}^3$ is obtained as an identifying space of a regular heptagon. 
It is triangulated as in the right-hand side of Figure \ref{Figure7}. 

\begin{figure}[hbtp]
\setlength{\unitlength}{1cm}
\begin{center}
\includegraphics[height=6cm]{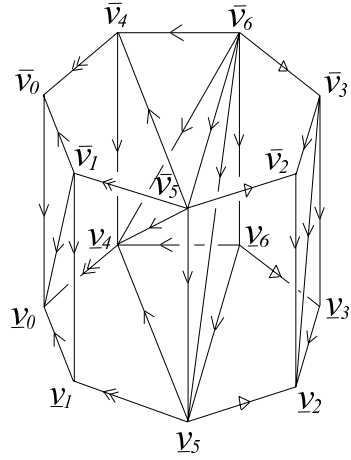}
\caption{a triangulation of $Y\times S^1$ \label{Figure8}}
\end{center}
\end{figure}

Therefore, the $3$-manifold $Y\times S^1$ is realized in $\mathbb{R}^3$ 
by the (singular) locally ordered complex $\tilde{\mathcal{K}}$ illustrated in Figure \ref{Figure8} 
(the cube with vertices $\underline{v}_0,\underline{v}_1,\underline{v}_4,\underline{v}_5$, 
$\overline{v}_0,\overline{v}_1,\overline{v}_4,\overline{v}_5$ and 
the cube with vertices $\underline{v}_2,\underline{v}_3,\underline{v}_5,\underline{v}_6$, 
$\overline{v}_2,\overline{v}_3,\overline{v}_5,\overline{v}_6$ are decomposed as product complexes, respectively). 

Hence, the product of the fusion algebra of $Z^{\Delta }$ is given by the linear map 
$$Z_{\text{\scriptsize fusion}}^{\Delta }:Z^{\Delta }(S^1\times S^1)\otimes Z^{\Delta }(S^1\times S^1)
\longrightarrow Z^{\Delta }(S^1\times S^1),$$
which is induced from the cobordism with triangulated boundary
$$W_{\text{\scriptsize fusion}}:=(Y\times S^1;(C_1\times S^1,\mathcal{K}_1)\cup (C_2\times S^1,\mathcal{K}_2), (C_3\times S^1,\mathcal{K}_3)),$$
where $\mathcal{K}_i\ (i=1,2,3)$ is the sub-complex of $\tilde{\mathcal{K}}$ such that 
\begin{align*}
\mathcal{K}_1&=K(|\underline{v_0}\ \overline{v_0}\ \overline{v_1}|)\cup K(|\underline{v_0}\ \underline{v_1}\ \overline{v_1}|) ,\\ 
\mathcal{K}_2&=K(|\underline{v_2}\ \overline{v_2}\ \overline{v_3}|)\cup K(|\underline{v_2}\ \underline{v_3}\ \overline{v_3}|), \\ 
\mathcal{K}_3&=K(|\underline{v_4}\ \overline{v_4}\ \overline{v_6}|)\cup K(|\underline{v_4}\ \underline{v_6}\ \overline{v_6}|) . \end{align*}

We note that $\mathcal{K}_i$ is a copy of $\mathcal{K}$ for each $i=1,2,3$.
\par
\begin{figure}[hbtp]
\setlength{\unitlength}{1cm}
\begin{center}
\includegraphics[height=5cm]{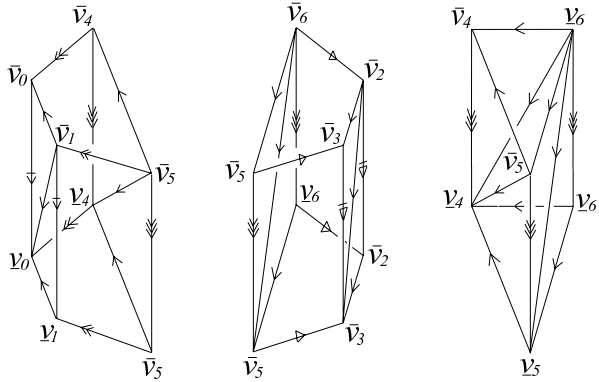}
\caption{three sub-complexes of $\tilde{\mathcal{K}}$ \label{Figure9}}
\end{center}
\end{figure}

\begin{figure}[hbt]
\setlength{\unitlength}{1cm}
\begin{center}
\includegraphics[height=4.5cm]{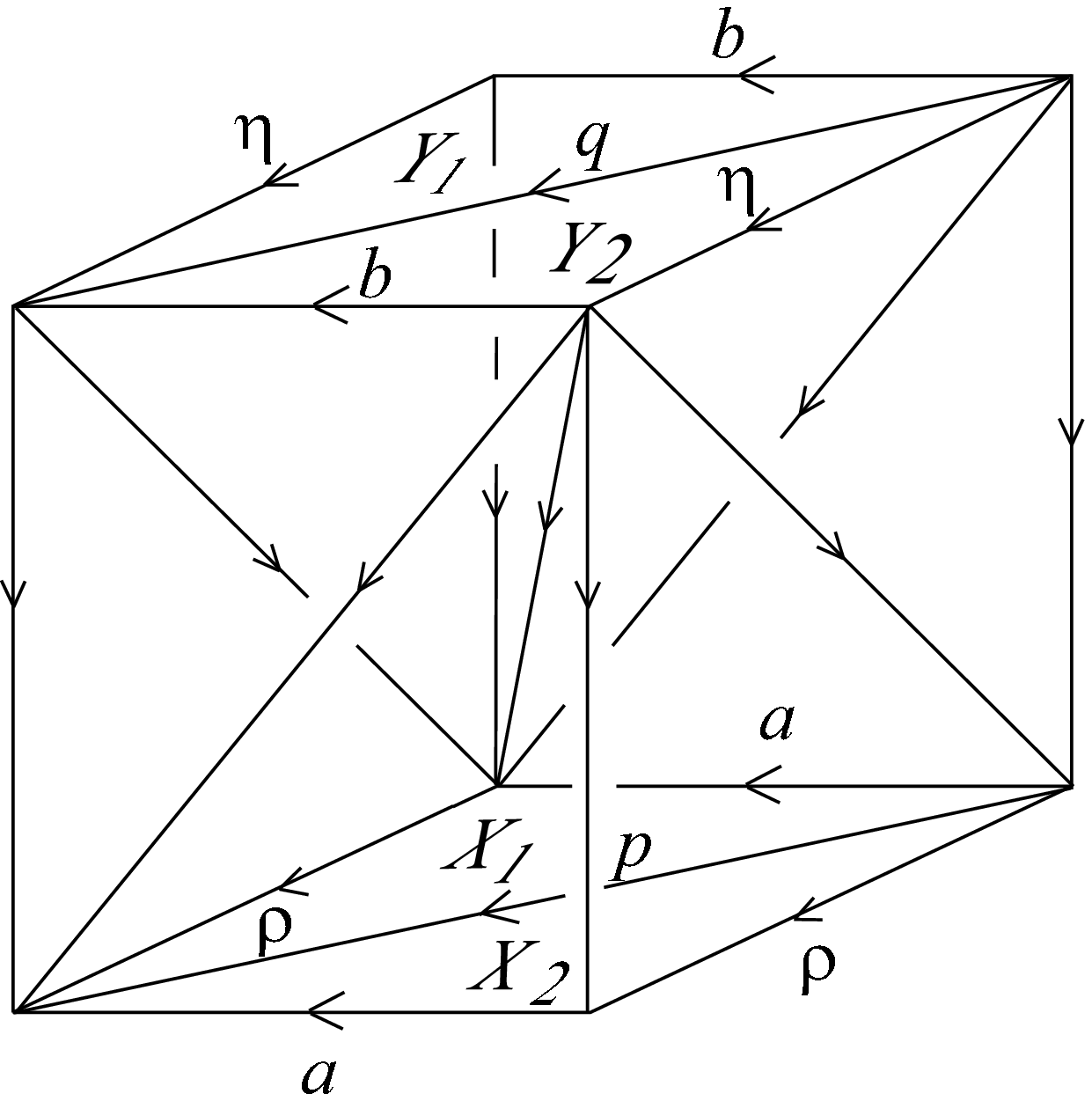}
\qquad \quad 
\includegraphics[height=4.5cm]{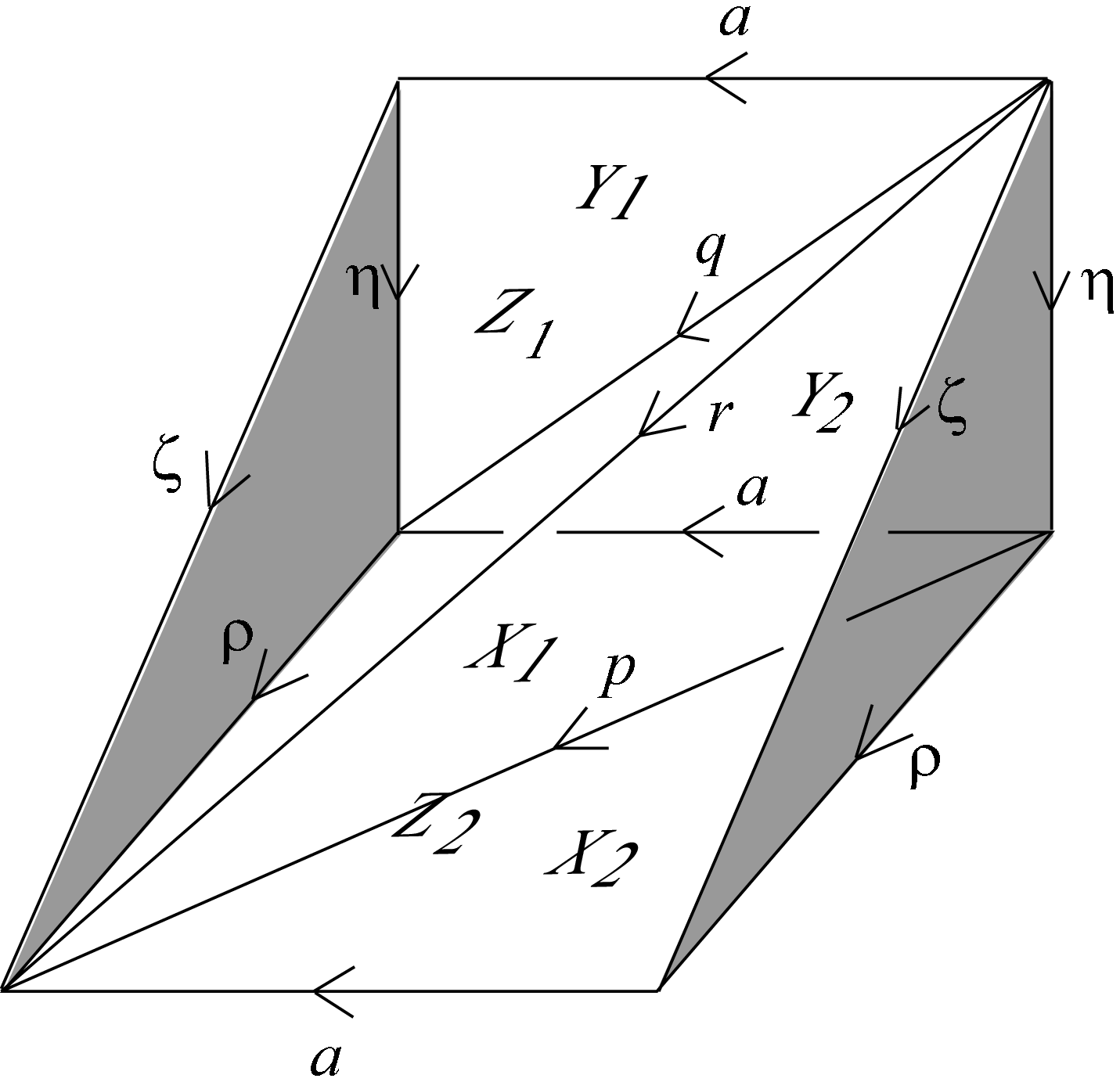}
\caption{triangulations of $S^1\times S^1\times [0,1]$ and $H\times S^1$ with colored boundaries \label{Figure10}}
\end{center}
\end{figure}

 
Let $\tilde{\mathcal{K}}_1, \tilde{\mathcal{K}}_2, \mathcal{L}$ be the sub-complexes of $\tilde{\mathcal{K}}$ 
depicted in Figure \ref{Figure9}. 
The geometrical realizations of $\tilde{\mathcal{K}}_1$ and $\tilde{\mathcal{K}}_2$ are $S^1\times S^1\times [0,1]$, 
and the geometrical realization of $\mathcal{L}$ is $H\times S^1$, 
where $H$ is the quotient space obtained from a $2$-simplex by identifying its three vertices. 
Let $Z(S^1\times S^1\times [0,1];{}_{\rho ,a,p,X_1,X_2}^{\eta ,b,q,Y_1,Y_2})$ and 
$Z(H\times S^1;a,{}_{\rho , p,X_1,X_2;\eta ,q,Y_1,Y_2}^{\zeta ,r,Z_1,Z_2})$ be the Turaev-Viro-Ocneanu invariants 
of $S^1\times S^1\times [0,1]$ triangulated by $\mathcal{K}\times [0,1]$ and 
of $H\times S^1$ triangulated by $\mathcal{L}$ assigned colors to their boundaries 
as in Figure \ref{Figure10}, respectively. 

\begin{figure}[hbtp]
\setlength{\unitlength}{0.7cm}
\begin{center}
\begin{picture}(17,3.5)
\put(1.0,1){\framebox(2,2)}
\put(1.0,1){\line(1,1){2}}
\put(1.9,1.9){\line(1,0){0.2}}
\put(1.9,1.9){\line(0,1){0.2}}
\put(2.1,2.4){$p$}
\put(1.8,0.4){$a$}
\put(1.8,3.2){$a$}
\put(0.45,2.0){$\rho $}
\put(3.2,2.0){$\rho $}
\put(1.7,0.85){$\ll $}
\put(1.7,2.86){$\ll $}
\put(0.8,2){$\vee $}
\put(2.82,2){$\vee $}
\put(1.3,2.4){$X_1$}
\put(2.1,1.4){$X_2$}
\put(1.8,-0.3){$\mathcal{K}_1$}
\put(6.5,1){\framebox(2,2)}
\put(6.5,1){\line(1,1){2}}
\put(7.4,1.9){\line(1,0){0.2}}
\put(7.4,1.9){\line(0,1){0.2}}
\put(7.6,2.4){$q$}
\put(7.3,0.4){$b$}
\put(7.3,3.2){$b$}
\put(5.95,2.0){$\eta $}
\put(8.7,2.0){$\eta $}
\put(7.2,0.85){$\ll $}
\put(7.2,2.86){$\ll $}
\put(6.3,2){$\vee $}
\put(8.32,2){$\vee $}
\put(6.8,2.4){$Y_1$}
\put(7.6,1.4){$Y_2$}
\put(7.3,-0.3){$\mathcal{K}_2$}
\put(12.2,1){\framebox(2,2)}
\put(12.2,1){\line(1,1){2}}
\put(13.1,1.9){\line(1,0){0.2}}
\put(13.1,1.9){\line(0,1){0.2}}
\put(13.3,2.4){$r$}
\put(13.0,0.4){$c$}
\put(13.0,3.4){$c$}
\put(11.65,2.0){$\zeta $}
\put(14.4,2.0){$\zeta $}
\put(12.9,0.85){$\ll $}
\put(12.9,2.86){$\ll $}
\put(12.0,2){$\vee $}
\put(14.02,2){$\vee $}
\put(12.5,2.4){$Z_1$}
\put(13.3,1.4){$Z_2$}
\put(13.0,-0.3){$\mathcal{K}_3$}
\end{picture}
\caption{colors of $\mathcal{K}_1$, $\mathcal{K}_2$ and $\mathcal{K}_3$ \label{Figure11}}
\end{center}
\end{figure}
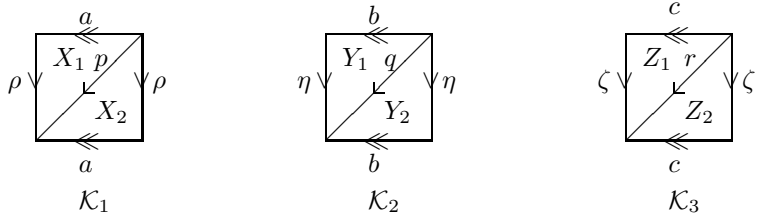

\par 
\pagebreak 

Then, the Turaev-Viro-Ocneanu invariant of $Y\times S^1$ triangulated by 
$\tilde{\mathcal{K}}$ assigned colors as in Figure \ref{Figure11} to 
$\mathcal{K}_1,\ \mathcal{K}_2,\ \mathcal{K}_3$ is given by 

\begin{align*}
&Z(Y\times S^1;{}_{\rho ,a,p,X_1,X_2;\eta ,b,q,Y_1,Y_2}^{\zeta ,c,r,Z_1,Z_2} )=
\dfrac{\lambda }{d(c)}\kern-0.9em \sum\limits _{\substack{\rho ',\eta ', p',q'\in \Delta \\ X_1'\in \mathcal{B}_{\rho 'c}^{p'}, X_2'\in \mathcal{B}_{c\rho '}^{p'}\\ Y_1'\in \mathcal{B}_{\eta 'c}^{q'}, Y_2'\in \mathcal{B}_{c\eta '}^{q'}}}
\kern-0.9em Z(H\times S^1;c,{}_{\rho ',p',X_1',X_2';\eta ',q',Y_1',Y_2'}^{\zeta ,r,Z_1,Z_2}) \notag \\ 
&\qquad \qquad \quad \qquad \times Z(S^1\times S^1\times [0,1];{}_{\rho ,a,p,X_1,X_2}^{\rho ',c,p',X_1',X_2'}) 
Z(S^1\times S^1\times [0,1];{}_{\eta ,b,q,Y_1,Y_2}^{\eta ',c,q',Y_1',Y_2'}). 
\end{align*}

Since 

\par 
\begin{figure}[hbtp]
\setlength{\unitlength}{0.7cm}
\begin{center}
\begin{picture}(17,3.2)
\put(1.5,1){\framebox(2,2)}
\put(1.5,1){\line(1,1){2}}
\put(2.4,1.9){\line(1,0){0.2}}
\put(2.4,1.9){\line(0,1){0.2}}
\put(2.6,2.4){$p$}
\put(2.3,0.4){$a$}
\put(2.3,3.2){$a$}
\put(0.95,2.0){$\rho $}
\put(3.7,2.0){$\rho $}
\put(2.2,0.85){$\ll $}
\put(2.2,2.86){$\ll $}
\put(1.3,2){$\vee $}
\put(3.32,2){$\vee $}
\put(1.8,2.4){$X_1$}
\put(2.6,1.4){$X_2$}
\put(-1,1.9){$\Phi _{Id_{\mathcal{K},\mathcal{K}}}$}
\put(0.6,1.9){$\Biggl ($}
\put(4,1.9){$\Biggr ) $}
\put(4.65,1.9){$=$}
\put(4.8,1.9){$\sum\limits_{\substack{\eta ,b,q\in \Delta \\ Y_1\in \mathcal{B}_{\eta b}^q\\ Y_2\in \mathcal{B}_{b \eta }^q}} Z(S^1\times S^1\times [0,1];{}_{\rho ,a,p,X_1,X_2}^{\eta ,b,q,Y_1,Y_2})$}
\put(14.5,1){\framebox(2,2)}
\put(14.5,1){\line(1,1){2}}
\put(15.4,1.9){\line(1,0){0.2}}
\put(15.4,1.9){\line(0,1){0.2}}
\put(15.6,2.4){$q$}
\put(15.3,0.4){$b$}
\put(15.3,3.2){$b$}
\put(13.95,2.0){$\eta $}
\put(16.7,2.0){$\eta $}
\put(15.2,0.85){$\ll $}
\put(15.2,2.86){$\ll $}
\put(14.3,2){$\vee $}
\put(16.32,2){$\vee $}
\put(14.8,2.4){$Y_1$}
\put(15.6,1.4){$Y_2$}
\put(17,1.4){,}
\end{picture}
\end{center}
\end{figure}

\par \noindent 
and 
$$Z(H\times S^1;c,{}_{\rho ',p',X_1',X_2';\eta ',q',Y_1',Y_2'}^{\zeta ,r,Z_1,Z_2})
=\dfrac{1}{d(c)}\overline{Z(D^2\times S^1;{}_{c,\rho ',p',X_2',X_1';c,\eta ',q',Y_2',Y_1'}^{c,\zeta ,r,Z_2',Z_1'})},$$ 
it follows that the following diagram commutes. 

{\small $$\setlength{\unitlength}{0.7mm}
\begin{picture}(150,50)(0,0)
\put(0,0){\makebox(20,10)[c]{$V^{\Delta }(S^1\times S^1)\otimes V^{\Delta }(S^1\times S^1)$}}
\put(0,30){\makebox(20,10)[c]{$V^{\Delta }(S^1\times S^1)\otimes V^{\Delta }(S^1\times S^1)$}}
\put(135,30){\makebox(20,10)[l]{$V^{\Delta }(S^1\times S^1)$}}
\put(78,0){\makebox(20,10)[c]{$V^{\Delta }(S^1\times S^1)\otimes V^{\Delta }(S^1\times S^1)$}}
\put(135,0){\makebox(20,10)[l]{$V^{\Delta }(S^1\times S^1)$}}
\put(75,40){$\Phi _{W_{\text{fusion}}}$}
\put(-25,20){$\Phi _{Id_{\mathcal{K},\mathcal{K}}}\otimes \Phi _{Id_{\mathcal{K},\mathcal{K}}}$}
\put(153,20){$P$}
\put(45,35){\vector(1,0){85}}
\put(10,30){\vector(0,-1){20}}
\put(150,30){\vector(0,-1){20}}
\put(42,5){\vector(1,0){15}}
\put(119,5){\vector(1,0){15}}
\put(43,8){$P\otimes P$}
\put(123,8){$\tilde{m}$}
\end{picture}
$$}

This implies that the diagram in the proposition is commutative. 
This completes the proof. 
\qed 

\par \bigskip 
Let $S$ and $T$ be the orientation preserving diffeomorphisms on $S^1\times S^1$ as in (\ref{eq1}). 
We will compute the actions $Z^{\Delta }(S)$ and $Z^{\Delta }(T)$ on $Z^{\Delta }(S^1\times S^1)$ 
by using the singular triangulation $\mathcal{K}$ of $S^1\times S^1$ depicted in Figure \ref{Figure5}. 

\par \bigskip 
\begin{lem}
\label{Lemma2}  
Let $\Delta $ be a finite system of $\text{End}(M)_0$ obtained from a subfactor $N\subset M$ 
of an infinite factor $M$ with finite index and finite depth. 
Let $\tilde{S}_{\#}:V^{\Delta }(S^1\times S^1)\longrightarrow V^{\Delta }(S^1\times S^1)$ and 
$\tilde{T}_{\#}^{-1}:V^{\Delta }(S^1\times S^1)\longrightarrow V^{\Delta }(S^1\times S^1)$ be 
the $\mathbb{C}$-linear maps defined by 

\begin{align*}
\tilde{S}_{\#}(X_1\otimes X_2)
&=\sum\limits_{\substack{r\in \Delta \\ Z_1\in \mathcal{B}_{\bar{a}\rho }^r\\ Z_2\in \mathcal{B}_{\rho \bar{a}}^r}} 
\kern-0.2cm \dfrac{\sqrt{d(p)d(r)}}{d(\rho )}Z_2^{\ast }\tilde{X}_2^{\ast }\bar{a}(\hat{X}_1)Z_1\ Z_1\otimes Z_2 ,\\ 
\tilde{T}_{\#}^{-1}(X_1\otimes X_2)
&=\sum\limits_{\substack{r\in \Delta \\ Z_1\in \mathcal{B}_{\rho p}^r\\ Z_2\in \mathcal{B}_{p\rho }^r}}
\kern-0.2cm \dfrac{\sqrt{d(a)d(r)}}{d(p)}Z_2^{\ast }X_1^{\ast }\rho(X_2)Z_1\ Z_1\otimes Z_2 
\end{align*}
for all $X_1\in \mathcal{B}_{\rho a}^p$ and $X_2\in \mathcal{B}_{a\rho }^p$. 
Then, the actions $Z^{\Delta }(S)$ and $Z^{\Delta }(T)$ on $Z^{\Delta }(S^1\times S^1)$ 
are determined by the following commutative diagrams. 

\begin{equation}
\begin{CD}
V^{\Delta }(S^1\times S^1) @>\tilde{S}_{\#}>> V^{\Delta }(S^1\times S^1) \\ 
@VuVV @VVuV \\ 
Z^{\Delta }(S^1\times S^1) @>Z^{\Delta }(S)>> Z^{\Delta }(S^1\times S^1),
\end{CD}
\label{eq2} \end{equation}
\begin{equation}
\begin{CD}
V^{\Delta }(S^1\times S^1) @>\tilde{T}_{\#}^{-1}>> V^{\Delta }(S^1\times S^1) \\ 
@VuVV @VVuV \\ 
Z^{\Delta }(S^1\times S^1) @>Z^{\Delta }(T)^{-1}>> Z^{\Delta }(S^1\times S^1)
\end{CD}
\label{eq3} \end{equation}
Here, $u$ is the universal arrow associated to $Z^{\Delta }(S^1\times S^1)$ defined as in Lemma \ref{Lemma1}. 
\end{lem}

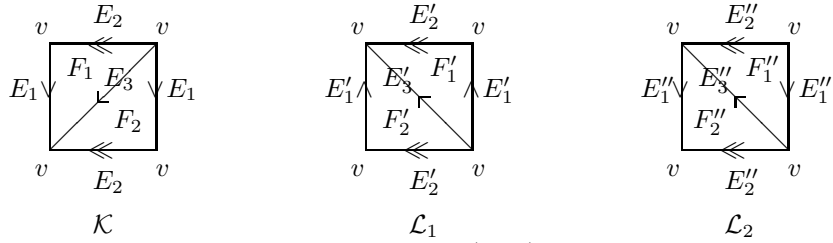
\begin{figure}[htbp]
\begin{center}
\setlength{\unitlength}{0.7cm}
\begin{picture}(14,3.7)
\put(1,1){\framebox(2,2)}
\put(1,1){\line(1,1){2}}
\put(1.9,1.9){\line(1,0){0.2}}
\put(1.9,1.9){\line(0,1){0.2}}
\put(0.7,0.5){$v$}
\put(3.0,0.5){$v$}
\put(0.7,3.2){$v$}
\put(3.0,3.2){$v$}
\put(2.0,2.2){$E_3$}
\put(1.8,0.25){$E_2$}
\put(1.8,3.4){$E_2$}
\put(0.2,2.0){$E_1$}
\put(3.2,2.0){$E_1$}
\put(1.7,0.85){$\ll $}
\put(1.7,2.85){$\ll $}
\put(0.79,2){$\vee $}
\put(2.82,2){$\vee $}
\put(1.3,2.4){$F_1$}
\put(2.2,1.4){$F_2$}
\put(1.8,-0.6){$\mathcal{K}$}
\put(7,1){\framebox(2,2)}
\put(7,3){\line(1,-1){2}}
\put(8,2.0){\line(1,0){0.2}}
\put(8,2.0){\line(0,-1){0.2}}
\put(6.7,0.5){$v$}
\put(9.0,0.5){$v$}
\put(6.7,3.2){$v$}
\put(9.0,3.2){$v$}
\put(7.3,2.2){$E_3'$}
\put(7.8,0.25){$E_2'$}
\put(7.8,3.4){$E_2'$}
\put(6.2,2.0){$E_1'$}
\put(9.2,2.0){$E_1'$}
\put(7.7,0.85){$\ll $}
\put(7.7,2.85){$\ll $}
\put(6.79,2){$\wedge $}
\put(8.82,2){$\wedge $}
\put(8.2,2.4){$F_1'$}
\put(7.3,1.4){$F_2'$}
\put(7.8,-0.6){$\mathcal{L}_1$}
\put(13,1){\framebox(2,2)}
\put(13,3){\line(1,-1){2}}
\put(14,2.0){\line(1,0){0.2}}
\put(14,2.0){\line(0,-1){0.2}}
\put(12.7,0.5){$v$}
\put(15.0,0.5){$v$}
\put(12.7,3.2){$v$}
\put(15.0,3.2){$v$}
\put(13.3,2.2){$E_3''$}
\put(13.8,0.25){$E_2''$}
\put(13.8,3.4){$E_2''$}
\put(12.2,2.0){$E_1''$}
\put(15.2,2.0){$E_1''$}
\put(13.7,0.85){$\ll $}
\put(13.7,2.85){$\ll $}
\put(12.79,2){$\vee $}
\put(14.82,2){$\vee $}
\put(14.2,2.4){$F_1''$}
\put(13.2,1.4){$F_2''$}
\put(13.8,-0.6){$\mathcal{L}_2$}
\end{picture}
\caption{three triangulations of $S^1\times S^1$ \label{Figure12}}
\end{center}
\end{figure}

\indent 
{\bf Proof.} 
Let $\mathcal{K}$, $\mathcal{L}_1$ and $\mathcal{L}_2$ be the locally ordered complexes depicted 
in Figure \ref{Figure12} that give singular triangulations of $S^1\times S^1$ (the opposite sides are identified).  
By considering the lifts of $S$ and $T$ to the universal covering
$$\mathbb{R}^2\longrightarrow S^1\times S^1,\ (x,y)\longmapsto (\text{exp}(2\pi iy),\text{exp}(2\pi ix)),$$
we see that $S$ and $T^{-1}$ are simplicial maps from $\mathcal{K}$ to $\mathcal{L}_1$ 
and from $\mathcal{K}$ to $\mathcal{L}_2$, respectively.  

\par \ 
\newpage 
So, for the color depicted in Figure \ref{Figure6}, we have

\begin{figure}[hbtp]
\setlength{\unitlength}{0.7cm}
\begin{center}
\begin{picture}(13,3.5)
\put(1.5,1){\framebox(2,2)}
\put(1.5,1){\line(1,1){2}}
\put(2.4,1.9){\line(1,0){0.2}}
\put(2.4,1.9){\line(0,1){0.2}}
\put(2.6,2.4){$p$}
\put(2.3,0.4){$a$}
\put(2.3,3.2){$a$}
\put(0.95,2.0){$\rho $}
\put(3.6,2.0){$\rho $}
\put(2.2,0.85){$\ll $}
\put(2.2,2.86){$\ll $}
\put(1.3,2){$\vee $}
\put(3.32,2){$\vee $}
\put(1.8,2.4){$X_1$}
\put(2.6,1.4){$X_2$}
\put(-0.9,1.9){$\Phi (S)$}
\put(0.6,1.9){$\Biggl ($}
\put(4,1.9){$\Biggr ) $}
\put(4.75,1.9){$=$}
\put(6,1){\framebox(2,2)}
\put(6,3){\line(1,-1){2}}
\put(7.0,2){\line(1,0){0.2}}
\put(7.0,2){\line(0,-1){0.2}}
\put(6.8,2.4){$p$}
\put(6.8,0.45){$\rho $}
\put(6.8,3.25){$\rho $}
\put(5.45,2.0){$a$}
\put(8.25,2.0){$a$}
\put(6.7,0.85){$<$}
\put(6.7,2.86){$<$}
\put(5.8,2){$\wedge $}
\put(5.8,2.2){$\wedge $}
\put(7.82,2){$\wedge $}
\put(7.82,2.2){$\wedge $}
\put(7.1,2.4){$X_1$}
\put(6.3,1.4){$X_2$}
\put(8.3,0.9){,}
\end{picture}
\end{center}
\begin{center}
\begin{picture}(13,3)
\put(1.5,1){\framebox(2,2)}
\put(1.5,1){\line(1,1){2}}
\put(2.4,1.9){\line(1,0){0.2}}
\put(2.4,1.9){\line(0,1){0.2}}
\put(2.6,2.4){$p$}
\put(2.3,0.4){$a$}
\put(2.3,3.2){$a$}
\put(0.95,2.0){$\rho $}
\put(3.6,2.0){$\rho $}
\put(2.2,0.85){$\ll $}
\put(2.2,2.86){$\ll $}
\put(1.3,2){$\vee $}
\put(3.32,2){$\vee $}
\put(1.8,2.4){$X_1$}
\put(2.6,1.4){$X_2$}
\put(-1.4,1.9){$\Phi (T^{-1})$}
\put(0.6,1.9){$\Biggl ($}
\put(4,1.9){$\Biggr ) $}
\put(4.75,1.9){$=$}
\put(6,1){\framebox(2,2)}
\put(6,3){\line(1,-1){2}}
\put(7.0,2){\line(1,0){0.2}}
\put(7.0,2){\line(0,-1){0.2}}
\put(7.08,1.92){\line(1,0){0.2}}
\put(7.08,1.92){\line(0,-1){0.2}}
\put(6.77,2.23){$a$}
\put(6.8,0.45){$p$}
\put(6.8,3.25){$p$}
\put(5.45,2.0){$\rho $}
\put(8.2,2.0){$\rho $}
\put(6.7,0.85){$<$}
\put(6.7,2.87){$<$}
\put(5.8,2){$\vee$}
\put(7.82,2){$\vee$}
\put(7.1,2.4){$X_2$}
\put(6.3,1.4){$X_1$}
\put(8.3,0.9){.}
\end{picture}
\end{center}
\end{figure}

\par 
Let $S_{\# }$ and $T_{\#}^{-1}$ be the linear maps from $V(S^1\times S^1;\mathcal{K})$ 
to $V(S^1\times S^1;\mathcal{K})$ defined by 
$$S_{\#}=\Phi _{\mathcal{L}_1,\mathcal{K}}\circ \Phi (S),\qquad T_{\#}^{-1}
=\Phi _{\mathcal{L}_2,\mathcal{K}}\circ \Phi (T^{-1}).$$
Then we have 

\begin{figure}[btp]
\vspace{-1.5cm}
\setlength{\unitlength}{0.7cm}
\begin{center}
\begin{picture}(13,3)
\put(1.5,1){\framebox(2,2)}
\put(1.5,1){\line(1,1){2}}
\put(2.4,1.9){\line(1,0){0.2}}
\put(2.4,1.9){\line(0,1){0.2}}
\put(2.6,2.4){$p$}
\put(2.3,0.4){$a$}
\put(2.3,3.2){$a$}
\put(0.95,2.0){$\rho $}
\put(3.7,2.0){$\rho $}
\put(2.2,0.85){$\ll $}
\put(2.2,2.86){$\ll $}
\put(1.3,2){$\vee $}
\put(3.32,2){$\vee $}
\put(1.8,2.4){$X_1$}
\put(2.6,1.4){$X_2$}
\put(-0.5,1.9){$S_{\#}$}
\put(0.6,1.9){$\Biggl ($}
\put(4,1.9){$\Biggr ) $}
\put(4.75,1.9){$=$}
\put(5.5,1.9){$\sum\limits_{\substack{\eta ,b,q\in \Delta \\ Y_1\in \mathcal{B}_{\eta b}^q\\ Y_2\in \mathcal{B}_{b \eta }^q}} Z(\mathcal{T}_1;{}_{\rho ,a,p,X_1,X_2}^{\eta ,b,q,Y_1,Y_2})$}
\put(12,1){\framebox(2,2)}
\put(12,1){\line(1,1){2}}
\put(12.9,1.9){\line(1,0){0.2}}
\put(12.9,1.9){\line(0,1){0.2}}
\put(13.1,2.4){$q$}
\put(12.8,0.4){$b$}
\put(12.8,3.2){$b$}
\put(11.45,2.0){$\eta $}
\put(14.2,2.0){$\eta $}
\put(12.7,0.85){$\ll $}
\put(12.7,2.86){$\ll $}
\put(11.8,2){$\vee $}
\put(13.82,2){$\vee $}
\put(12.3,2.4){$Y_1$}
\put(13.1,1.4){$Y_2$}
\put(14.5,1){,}
\end{picture}
\end{center}
\begin{center}
\begin{picture}(13,3)
\put(1.5,1){\framebox(2,2)}
\put(1.5,1){\line(1,1){2}}
\put(2.4,1.9){\line(1,0){0.2}}
\put(2.4,1.9){\line(0,1){0.2}}
\put(2.6,2.4){$p$}
\put(2.3,0.4){$a$}
\put(2.3,3.2){$a$}
\put(0.95,2.0){$\rho $}
\put(3.7,2.0){$\rho $}
\put(2.2,0.85){$\ll $}
\put(2.2,2.86){$\ll $}
\put(1.3,2){$\vee $}
\put(3.32,2){$\vee $}
\put(1.8,2.4){$X_1$}
\put(2.6,1.4){$X_2$}
\put(-0.6,1.9){$T_{\#}^{-1}$}
\put(0.6,1.9){$\Biggl ($}
\put(4,1.9){$\Biggr ) $}
\put(4.75,1.9){$=$}
\put(5.5,1.9){$\sum\limits_{\substack{\eta ,b,q\in \Delta \\ Y_1\in \mathcal{B}_{\eta b}^q\\ Y_2\in \mathcal{B}_{b \eta }^q}} Z(\mathcal{T}_2;{}_{\rho ,a,p,X_1,X_2}^{\eta ,b,q,Y_1,Y_2})$}
\put(12,1){\framebox(2,2)}
\put(12,1){\line(1,1){2}}
\put(12.9,1.9){\line(1,0){0.2}}
\put(12.9,1.9){\line(0,1){0.2}}
\put(13.1,2.4){$q$}
\put(12.8,0.4){$b$}
\put(12.8,3.2){$b$}
\put(11.45,2.0){$\eta $}
\put(14.2,2.0){$\eta $}
\put(12.7,0.85){$\ll $}
\put(12.7,2.86){$\ll $}
\put(11.8,2){$\vee $}
\put(13.82,2){$\vee $}
\put(12.3,2.4){$Y_1$}
\put(13.1,1.4){$Y_2$}
\put(14.5,1){,}
\end{picture}
\end{center}
\end{figure}

\begin{figure}[hbt]
\begin{center}
\setlength{\unitlength}{1cm}
\includegraphics[height=5cm]{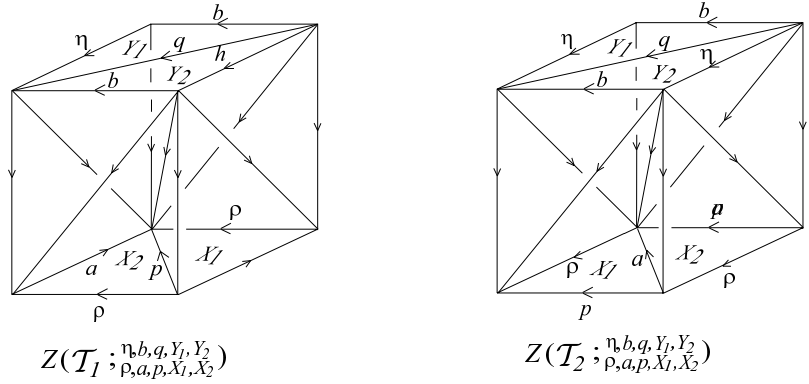}
\caption{two triangulations $\mathcal{T}_1$ and $\mathcal{T}_2$ of $S^1\times S^1\times [0,1]$ \label{Figure13}}
\end{center}
\end{figure}

\newpage 
\noindent 
where $Z(\mathcal{T}_1;{}_{\rho ,a,p,X_1,X_2}^{\eta ,b,q,Y_1,Y_2})$ is the Turaev-Viro-Ocneanu invariant 
of $(S^1\times S^1)\times [0,1]$ which boundary is triangulated by $\mathcal{K}$ and 
$\mathcal{L}_1$ assigned to its boundary the color depicted as in the left-hand side of Figure \ref{Figure13}, 
and $Z(\mathcal{T}_2;{}_{\rho ,a,p,X_1,X_2}^{\eta ,b,q,Y_1,Y_2})$ is the Turaev-Viro-Ocneanu invariant 
of $(S^1\times S^1)\times [0,1]$ which boundary is triangulated by $\mathcal{K}$ and 
$\mathcal{L}_2$ assigned to its boundary the color depicted as in the right-hand side of Figure \ref{Figure13}.

\par 
By tetrahedral symmetry on quantum $6j$-symbols, it can be proved that 
$$Z(\mathcal{T}_1;{}_{\rho ,a,p,X_1,X_2}^{\eta ,b,q,Y_1,Y_2})
=Z(\mathcal{T}_2;{}_{\bar{a} ,p,\rho ,\tilde{X}_2,\hat{X}_1}^{\eta ,b,q,Y_1,Y_2}).$$

\begin{figure}[hbt]
\begin{center}
\setlength{\unitlength}{1cm}
\includegraphics[height=4cm]{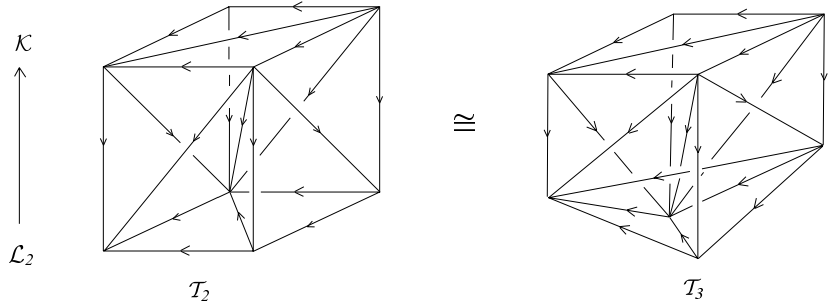}
\caption{an isomorphism between $(\mathcal{T}_2;\mathcal{L}_2,\mathcal{K})$ and 
$(\mathcal{T}_3;\mathcal{L}_2,\mathcal{K})$ \label{Figure14}}
\end{center}
\end{figure}

Let $(\mathcal{T}_3; \mathcal{L}_2,\mathcal{K})$ be the cobordism obtained by gluing one tetrahedron 
to the bottom of $\mathcal{K}\times [0,1]$. 
Then, we have 
$$Z(\mathcal{T}_2;{}_{\rho ,a,p,X_1,X_2}^{\eta ,b,q,Y_1,Y_2})
=Z(\mathcal{T}_3;{}_{\rho ,a,p,X_1,X_2}^{\eta ,b,q,Y_1,Y_2}),$$
since the cobordism $(\mathcal{T}_2;\mathcal{L}_2 ,\mathcal{K})$ 
is isomorphic to the cobordism $(\mathcal{T}_3; \mathcal{L}_2,\mathcal{K})$ (see Figure \ref{Figure14}). 

\par \newpage

Since 

\vspace{1.5cm}
\begin{figure}[ht]
\setlength{\unitlength}{1cm}
\begin{picture}(13,3.5)
\put(0,3.5){$Z(\mathcal{T}_3;{}_{\rho ,a,p,X_1,X_2}^{\eta ,b,q,Y_1,Y_2})
=\kern-1em\sum\limits_{\substack{r\in \Delta \\ Z_1\in \mathcal{B}_{\rho p}^r \\ Z_2\in \mathcal{B}_{p\rho }^r}}
\kern-0.5em\sqrt{d(a)d(r)}\ W\kern-0.3mm \Biggl ($}
\put(5.8,2.2){\includegraphics[width=2.5cm]{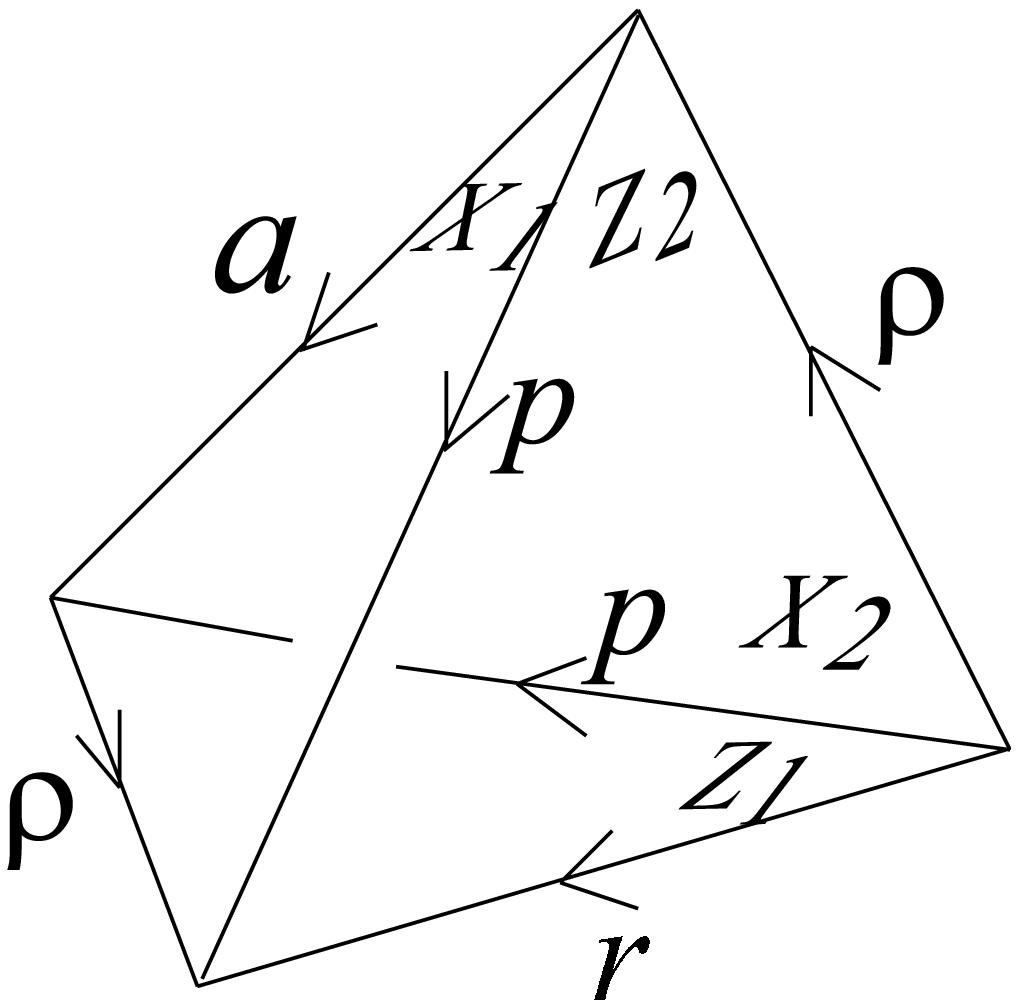}} 
\put(8.5,3.5){$\Biggr )\ Z\Biggl ($}
\put(9.4,2.0){\includegraphics[width=3cm]{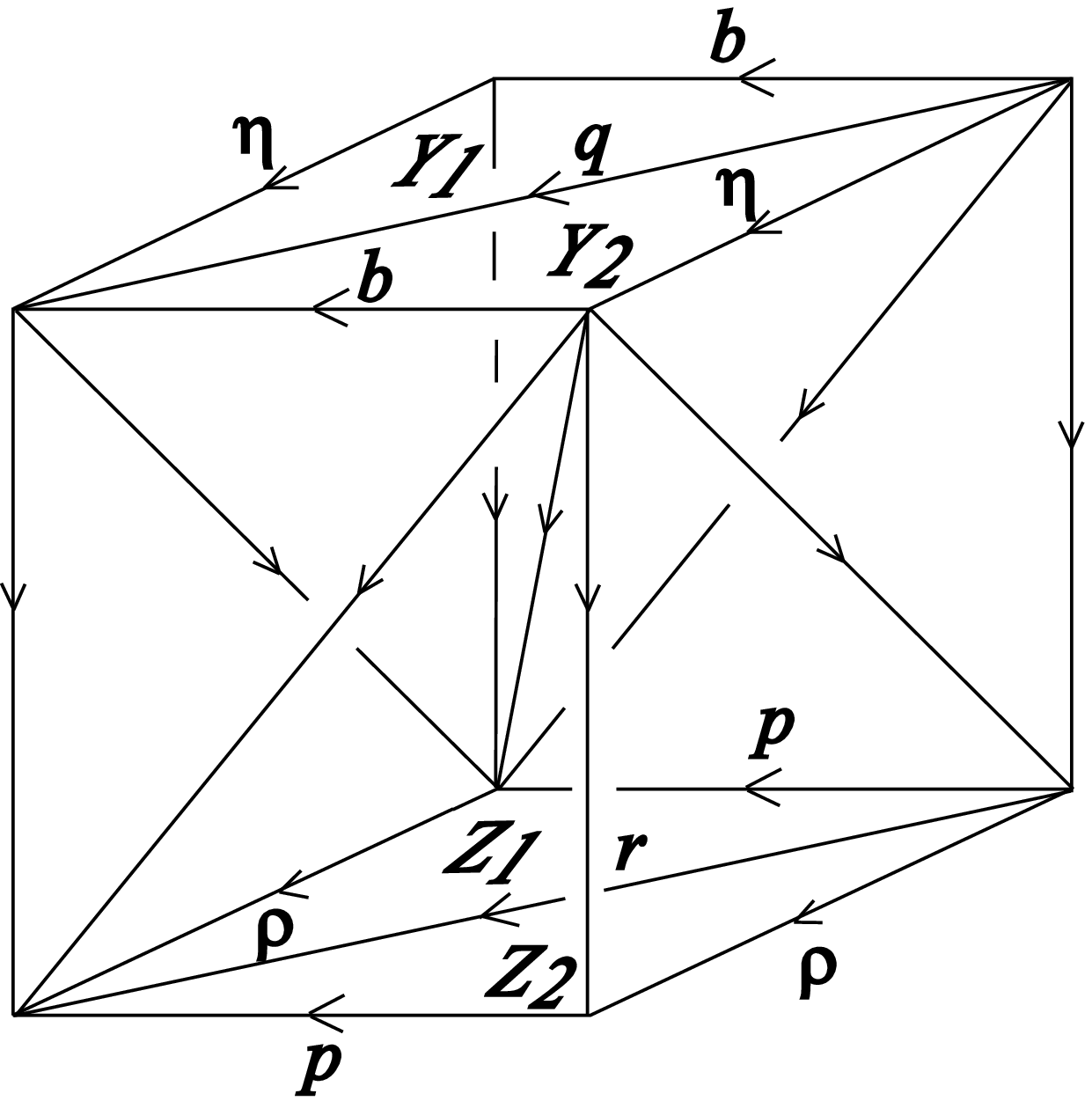}} 
\put(12.5,3.5){$\Biggr ),$}
\end{picture}
\end{figure}

\noindent 
it follows that 
\vspace{1cm}
\begin{figure}[htbp]
\setlength{\unitlength}{0.7cm}
\begin{center}
\begin{picture}(17,3)
\put(1.4,1){\framebox(2,2)}
\put(1.4,1){\line(1,1){2}}
\put(2.3,1.9){\line(1,0){0.2}}
\put(2.3,1.9){\line(0,1){0.2}}
\put(2.5,2.4){$p$}
\put(2.2,0.4){$a$}
\put(2.2,3.2){$a$}
\put(0.85,2.0){$\rho $}
\put(3.6,2.0){$\rho $}
\put(2.1,0.85){$\ll $}
\put(2.1,2.86){$\ll $}
\put(1.2,2){$\vee $}
\put(3.22,2){$\vee $}
\put(1.7,2.4){$X_1$}
\put(2.5,1.4){$X_2$}
\put(-0.6,1.9){$S_{\#}$}
\put(0.5,1.9){$\Biggl ($}
\put(4,1.9){$\Biggr ) $}
\put(4.75,1.9){$= \sum\limits_{\substack{r\in \Delta \\ Z_1\in \mathcal{B}_{\bar{a}\rho }^r\\ 
Z_2\in \mathcal{B}_{\rho \bar{a}}^r}} \dfrac{\sqrt{d(p)d(r)}}{d(\rho )}Z_2^{\ast }\tilde{X}_2^{\ast }\bar{a}
(\hat{X}_1)Z_1\ \Phi _{Id_{\mathcal{K},\mathcal{K}}}\kern0.1em \Biggl ( $}
\put(14.9,1){\framebox(2,2)}
\put(14.9,1){\line(1,1){2}}
\put(15.8,1.9){\line(1,0){0.2}}
\put(15.8,1.9){\line(0,1){0.2}}
\put(16.0,2.4){$r$}
\put(15.7,0.4){$\rho$}
\put(15.7,3.4){$\rho$}
\put(14.35,2.0){$\bar{a}$}
\put(17.1,2.0){$\bar{a}$}
\put(15.6,0.85){$\ll $}
\put(15.6,2.86){$\ll $}
\put(14.7,2){$\vee $}
\put(16.72,2){$\vee $}
\put(15.2,2.4){$Z_1$}
\put(16.0,1.4){$Z_2$}
\put(17.4,1.9){$\Biggr ) $}
\put(17.8,1.0){.}
\end{picture}
\end{center}
\end{figure}

\noindent 

This implies that $S_{\# }=\Phi _{Id_{\mathcal{K},\mathcal{K}}}\circ \tilde{S}_{\#}$, 
and whence, the diagram (\ref{eq2}) commutes. 
\par 
By a similar argument, it is also proved that 

\par \ 

\begin{figure}[hbt]
\setlength{\unitlength}{0.7cm}
\begin{center}
\begin{picture}(17,3)
\put(1.4,1){\framebox(2,2)}
\put(1.4,1){\line(1,1){2}}
\put(2.3,1.9){\line(1,0){0.2}}
\put(2.3,1.9){\line(0,1){0.2}}
\put(2.5,2.4){$p$}
\put(2.2,0.4){$a$}
\put(2.2,3.2){$a$}
\put(0.85,2.0){$\rho $}
\put(3.6,2.0){$\rho $}
\put(2.1,0.85){$\ll $}
\put(2.1,2.86){$\ll $}
\put(1.2,2){$\vee $}
\put(3.22,2){$\vee $}
\put(1.7,2.4){$X_1$}
\put(2.5,1.4){$X_2$}
\put(-0.6,1.9){$T^{-1}_{\#}$}
\put(0.5,1.9){$\Biggl ($}
\put(4,1.9){$\Biggr ) $}
\put(4.75,1.9){$=\kern-0.2em \sum\limits_{\substack{r\in \Delta \\ Z_1\in \mathcal{B}_{\rho p}^r\\ 
Z_2\in \mathcal{B}_{p\rho }^r}} \dfrac{\sqrt{d(a)d(r)}}{d(p)}Z_2^{\ast }X_1^{\ast }\rho (X_2)Z_1\ 
\Phi _{Id_{\mathcal{K},\mathcal{K}}}\Biggl ( $}
\put(14.7,1){\framebox(2,2)}
\put(14.7,1){\line(1,1){2}}
\put(15.6,1.9){\line(1,0){0.2}}
\put(15.6,1.9){\line(0,1){0.2}}
\put(15.8,2.4){$r$}
\put(15.5,0.4){$p$}
\put(15.5,3.4){$p$}
\put(14.15,2.0){$\rho $}
\put(16.9,2.0){$\rho $}
\put(15.4,0.85){$\ll $}
\put(15.4,2.86){$\ll $}
\put(14.5,2){$\vee $}
\put(16.52,2){$\vee $}
\put(15.0,2.4){$Z_1$}
\put(15.8,1.4){$Z_2$}
\put(17.3,1.9){$\Biggr ) $}
\put(17.7,1.0){,}
\end{picture}
\end{center}
\end{figure}

\noindent 
and whence, $T^{-1}_{\# }=\Phi _{Id_{\mathcal{K},\mathcal{K}}}\circ \tilde{T}_{\#}^{-1}$.  
This proves that the diagram (\ref{eq3}) commutes. 
\qed

\par \bigskip 
We note that the group $SL(2,\mathbb{Z})$ is not only generated by $S=\begin{pmatrix} 0 & 1 \\ -1 & 0
\end{pmatrix}$ and $T=\begin{pmatrix} 1 & 0 \\ 1 & 1\end{pmatrix}$, 
but also generated by $S'=\begin{pmatrix} 0 & -1 \\ 1 & 0
\end{pmatrix}$ and $T'=\begin{pmatrix} 1 & 1 \\ 0 & 1\end{pmatrix}$. 
The matrices $S'$ and $T'$ correspond to the orientation preserving diffeomorphisms 
from $S^1\times S^1$ to $S^1\times S^1$ defined by $S'(z,w)=(w,\bar{z})$ and $T'(z,w)=(z,zw)$ 
for all $(z,w)\in S^1\times S^1$, respectively. 
Since $S'$ and $T'$ has also the relations $(S')^4=I,\ (S'T')^3=(S')^2$, 
we have a group isomorphism $Q:SL(2,\mathbb{Z})\longrightarrow SL(2,\mathbb{Z})$ such that $Q(S)=S'$ and $Q(T)=T'$. 

\par \bigskip 
\begin{lem}
\label{Lemma3}  
Let $\Delta $ be a finite system of $\text{End}(M)_0$ obtained from a subfactor $N\subset M$ 
of an infinite factor $M$ with finite index and finite depth. 
Let $\bar{P}:Z^{\Delta }(S^1\times S^1)\longrightarrow Z^{\Delta }(S^1\times S^1)$ be 
the conjugate linear isomorphism induced from $P$ defined as in Proposition \ref{Proposition1}. 
Then, the following diagram commutes for all $R\in SL(2,\mathbb{Z})$. 
\begin{equation*}
\begin{CD}
Z^{\Delta }(S^1\times S^1) @>Z^{\Delta }(R)>> Z^{\Delta }(S^1\times S^1) \\ 
@V\bar{P}VV @VV\bar{P}V \\ 
Z^{\Delta }(S^1\times S^1) @>Z^{\Delta }(R')>> Z^{\Delta }(S^1\times S^1) 
\end{CD}
\end{equation*}
Here, $R'$ denotes the image of $R$ by the group isomorphism $Q:SL(2,\mathbb{Z})\longrightarrow SL(2,\mathbb{Z})$. 
\end{lem}

\indent 
{\bf Proof.} 
Let $\mathcal{K}$, $\mathcal{L}_1'$ and $\mathcal{L}_2'$ be the locally ordered complexes 
depicted in Figure \ref{Figure15} that give singular triangulations of $S^1\times S^1$ 
(the opposite sides are identified).  

\begin{figure}[hbtp]
\begin{center}
\setlength{\unitlength}{0.7cm}
\begin{picture}(14,3.7)
\put(1,1){\framebox(2,2)}
\put(1,1){\line(1,1){2}}
\put(1.9,1.9){\line(1,0){0.2}}
\put(1.9,1.9){\line(0,1){0.2}}
\put(0.7,0.5){$v$}
\put(3.0,0.5){$v$}
\put(0.7,3.2){$v$}
\put(3.0,3.2){$v$}
\put(2.0,2.2){$E_3$}
\put(1.8,0.25){$E_2$}
\put(1.8,3.4){$E_2$}
\put(0.2,2.0){$E_1$}
\put(3.2,2.0){$E_1$}
\put(1.7,0.85){$\ll $}
\put(1.7,2.85){$\ll $}
\put(0.79,2){$\vee $}
\put(2.82,2){$\vee $}
\put(1.3,2.4){$F_1$}
\put(2.2,1.4){$F_2$}
\put(1.8,-0.6){$\mathcal{K}$}
\put(7,1){\framebox(2,2)}
\put(7,3){\line(1,-1){2}}
\put(8,2.0){\line(-1,0){0.2}}
\put(8,2.0){\line(0,1){0.2}}
\put(6.7,0.5){$v$}
\put(9.0,0.5){$v$}
\put(6.7,3.2){$v$}
\put(9.0,3.2){$v$}
\put(8.0,1.5){$E_3'$}
\put(7.8,0.25){$E_2'$}
\put(7.8,3.4){$E_2'$}
\put(6.2,2.0){$E_1'$}
\put(9.2,2.0){$E_1'$}
\put(7.7,0.85){$\gg $}
\put(7.7,2.85){$\gg $}
\put(6.79,2){$\vee $}
\put(8.82,2){$\vee $}
\put(8.2,2.4){$F_1'$}
\put(7.3,1.4){$F_2'$}
\put(7.8,-0.6){$\mathcal{L}_1'$}
\put(13,1){\framebox(2,2)}
\put(13,3){\line(1,-1){2}}
\put(14,2.0){\line(-1,0){0.2}}
\put(14,2.0){\line(0,1){0.2}}
\put(12.7,0.5){$v$}
\put(15.0,0.5){$v$}
\put(12.7,3.2){$v$}
\put(15.0,3.2){$v$}
\put(14.0,1.5){$E_3''$}
\put(13.8,0.25){$E_2''$}
\put(13.8,3.4){$E_2''$}
\put(12.2,2.0){$E_1''$}
\put(15.2,2.0){$E_1''$}
\put(13.7,0.85){$\ll $}
\put(13.7,2.85){$\ll $}
\put(12.79,2){$\vee $}
\put(14.82,2){$\vee $}
\put(14.2,2.4){$F_1''$}
\put(13.2,1.4){$F_2''$}
\put(13.8,-0.6){$\mathcal{L}_2'$}
\end{picture}
\caption{three triangulations of $S^1\times S^1$ \label{Figure15}}
\end{center}
\end{figure}

\par 
Since $S'$ and ${T'}^{-1}$ are simplicial maps from $\mathcal{K}$ to $\mathcal{L}_1'$ 
and from $\mathcal{K}$ to $\mathcal{L}_2'$, respectively, we have 

\begin{figure}[hbtp]
\setlength{\unitlength}{0.7cm}
\begin{center}
\begin{picture}(13,3.5)
\put(1.5,1){\framebox(2,2)}
\put(1.5,1){\line(1,1){2}}
\put(2.4,1.9){\line(1,0){0.2}}
\put(2.4,1.9){\line(0,1){0.2}}
\put(2.6,2.4){$p$}
\put(2.3,0.4){$a$}
\put(2.3,3.2){$a$}
\put(0.95,2.0){$\rho $}
\put(3.6,2.0){$\rho $}
\put(2.2,0.85){$\ll $}
\put(2.2,2.86){$\ll $}
\put(1.3,2){$\vee $}
\put(3.32,2){$\vee $}
\put(1.8,2.4){$X_1$}
\put(2.6,1.4){$X_2$}
\put(-0.9,1.9){$\Phi (S')$}
\put(0.6,1.9){$\Biggl ($}
\put(4,1.9){$\Biggr ) $}
\put(4.75,1.9){$=$}
\put(6,1){\framebox(2,2)}
\put(6,3){\line(1,-1){2}}
\put(7.0,2){\line(-1,0){0.2}}
\put(7.0,2){\line(0,1){0.2}}
\put(6.8,2.4){$p$}
\put(6.8,0.45){$\rho $}
\put(6.8,3.25){$\rho $}
\put(5.45,2.0){$a$}
\put(8.25,2.0){$a$}
\put(6.7,0.85){$>$}
\put(6.7,2.86){$>$}
\put(5.8,2){$\vee $}
\put(5.8,2.2){$\vee $}
\put(7.82,2){$\vee $}
\put(7.82,2.2){$\vee $}
\put(7.1,2.4){$X_2$}
\put(6.3,1.4){$X_1$}
\put(8.3,0.9){,}
\end{picture}
\end{center}
\begin{center}
\begin{picture}(13,3)
\put(1.5,1){\framebox(2,2)}
\put(1.5,1){\line(1,1){2}}
\put(2.4,1.9){\line(1,0){0.2}}
\put(2.4,1.9){\line(0,1){0.2}}
\put(2.6,2.4){$p$}
\put(2.3,0.4){$a$}
\put(2.3,3.2){$a$}
\put(0.95,2.0){$\rho $}
\put(3.6,2.0){$\rho $}
\put(2.2,0.85){$\ll $}
\put(2.2,2.86){$\ll $}
\put(1.3,2){$\vee $}
\put(3.32,2){$\vee $}
\put(1.8,2.4){$X_1$}
\put(2.6,1.4){$X_2$}
\put(-1.5,1.9){$\Phi ({T'}^{-1})$}
\put(0.6,1.9){$\Biggl ($}
\put(4,1.9){$\Biggr ) $}
\put(4.75,1.9){$=$}
\put(6,1){\framebox(2,2)}
\put(6,3){\line(1,-1){2}}
\put(7.0,2){\line(-1,0){0.2}}
\put(7.0,2){\line(0,1){0.2}}
\put(6.77,2.23){$\rho $}
\put(6.8,0.45){$a$}
\put(6.8,3.25){$a$}
\put(5.45,2.0){$p$}
\put(8.2,2.0){$p$}
\put(6.7,0.85){$<$}
\put(6.7,2.87){$<$}
\put(5.8,2){$\vee$}
\put(7.82,2){$\vee$}
\put(7.1,2.4){$X_2$}
\put(6.3,1.4){$X_1$}
\put(8.3,0.9){.}
\end{picture}
\end{center}
\end{figure}

\par \newpage 
Thus, we have 

\begin{figure}[htbp]
\setlength{\unitlength}{0.7cm}
\begin{center}
\begin{picture}(13,3)
\put(-0.5,1){\framebox(2,2)}
\put(-0.5,1){\line(1,1){2}}
\put(0.4,1.9){\line(1,0){0.2}}
\put(0.4,1.9){\line(0,1){0.2}}
\put(0.6,2.4){$p$}
\put(0.3,0.4){$a$}
\put(0.3,3.2){$a$}
\put(-1.05,2.0){$\rho $}
\put(1.7,2.0){$\rho $}
\put(0.2,0.87){$\ll $}
\put(0.2,2.88){$\ll $}
\put(-0.68,2){$\vee $}
\put(1.34,2){$\vee $}
\put(-0.2,2.4){$X_1$}
\put(0.6,1.4){$X_2$}
\put(-2.5,1.9){$S_{\#}'$}
\put(-1.4,1.9){$\Biggl ($}
\put(2,1.9){$\Biggr ) $}
\put(2.75,1.9){$=$}
\put(3.5,1.9){$\sum\limits_{\substack{\eta ,b,q\in \Delta \\ Y_1\in \mathcal{B}_{\eta b}^q\\ 
Y_2\in \mathcal{B}_{b \eta }^q}} Z(\mathcal{T}_1';{}_{\rho ,a,p,X_1,X_2}^{\eta ,b,q,Y_1,Y_2})$}
\put(9.3,1.9){$\Phi _{Id_{\mathcal{K},\mathcal{K}}}$}
\put(10.75,1.9){$\Biggl ($}
\put(14.1,1.9){$\Biggr ) $}
\put(11.6,1){\framebox(2,2)}
\put(11.6,1){\line(1,1){2}}
\put(12.5,1.9){\line(1,0){0.2}}
\put(12.5,1.9){\line(0,1){0.2}}
\put(12.7,2.4){$q$}
\put(12.4,0.4){$b$}
\put(12.4,3.2){$b$}
\put(11.05,2.0){$\eta $}
\put(13.8,2.0){$\eta $}
\put(12.3,0.87){$\ll $}
\put(12.3,2.88){$\ll $}
\put(11.42,2){$\vee $}
\put(13.44,2){$\vee $}
\put(11.9,2.4){$Y_1$}
\put(12.7,1.4){$Y_2$}
\put(14.5,1){,}
\end{picture}
\end{center}
\begin{center}
\begin{picture}(13,3)
\put(-0.5,1){\framebox(2,2)}
\put(-0.5,1){\line(1,1){2}}
\put(0.4,1.9){\line(1,0){0.2}}
\put(0.4,1.9){\line(0,1){0.2}}
\put(0.6,2.4){$p$}
\put(0.3,0.4){$a$}
\put(0.3,3.2){$a$}
\put(-1.05,2.0){$\rho $}
\put(1.7,2.0){$\rho $}
\put(0.2,0.87){$\ll $}
\put(0.2,2.88){$\ll $}
\put(-0.68,2){$\vee $}
\put(1.34,2){$\vee $}
\put(-0.2,2.4){$X_1$}
\put(0.6,1.4){$X_2$}
\put(-2.6,1.9){${T_{\#}'}^{-1}$}
\put(-1.4,1.9){$\Biggl ($}
\put(2,1.9){$\Biggr ) $}
\put(2.75,1.9){$=$}
\put(3.5,1.9){$\sum\limits_{\substack{\eta ,b,q\in \Delta \\ 
Y_1\in \mathcal{B}_{\eta b}^q\\ Y_2\in \mathcal{B}_{b \eta }^q}} 
Z(\mathcal{T}_2';{}_{\rho ,a,p,X_1,X_2}^{\eta ,b,q,Y_1,Y_2})$}
\put(9.3,1.9){$\Phi _{Id_{\mathcal{K},\mathcal{K}}}$}
\put(10.75,1.9){$\Biggl ($}
\put(14.1,1.9){$\Biggr ) $}
\put(11.6,1){\framebox(2,2)}
\put(11.6,1){\line(1,1){2}}
\put(12.5,1.9){\line(1,0){0.2}}
\put(12.5,1.9){\line(0,1){0.2}}
\put(12.7,2.4){$q$}
\put(12.4,0.4){$b$}
\put(12.4,3.2){$b$}
\put(11.05,2.0){$\eta $}
\put(13.8,2.0){$\eta $}
\put(12.3,0.87){$\ll $}
\put(12.3,2.88){$\ll $}
\put(11.42,2){$\vee $}
\put(13.44,2){$\vee $}
\put(11.9,2.4){$Y_1$}
\put(12.7,1.4){$Y_2$}
\put(14.5,1){,}
\end{picture}
\end{center}
\end{figure}
\begin{figure}[hbt]
\begin{center}
\setlength{\unitlength}{1cm}
\includegraphics[height=5cm]{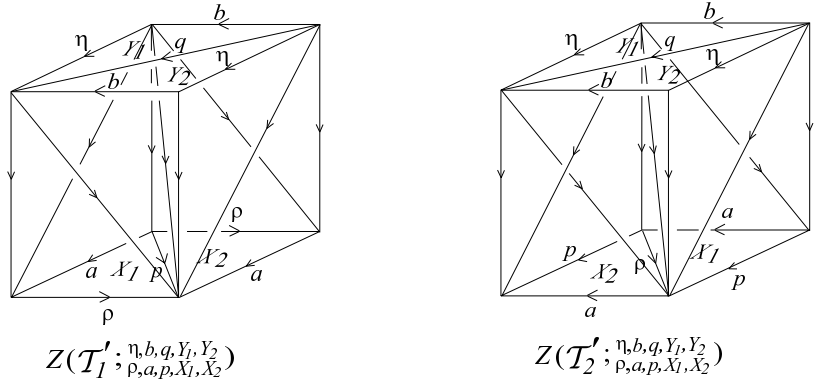}
\caption{two triangulations $\mathcal{T}_1'$ and $\mathcal{T}_2'$ of $S^1\times S^1\times [0,1]$ \label{Figure16}}
\end{center}
\end{figure}


\noindent 
where $Z(\mathcal{T}_1';{}_{\rho ,a,p,X_1,X_2}^{\eta ,b,q,Y_1,Y_2})$ 
and $Z(\mathcal{T}_2';{}_{\rho ,a,p,X_1,X_2}^{\eta ,b,q,Y_1,Y_2})$ are 
the Turaev-Viro-Ocneanu invariants of $S^1\times S^1\times [0,1]$ based on the triangulations $\mathcal{T}_1'$ and 
$\mathcal{T}_2'$ with colored boundaries as in Figure \ref{Figure16}, respectively. 

\par 
Since 
$$Z(\mathcal{T}_1';{}_{\rho ,a,p,X_1,X_2}^{\eta ,b,q,Y_1,Y_2})
=\overline{Z(\mathcal{T}_1;{}_{a,\rho ,p,X_2,X_1}^{b,\eta ,q,Y_2,Y_1})} 
\text{ and } Z(\mathcal{T}_2';{}_{\rho ,a,p,X_1,X_2}^{\eta ,b,q,Y_1,Y_2})
=\overline{Z(\mathcal{T}_2;{}_{a,\rho ,p,X_2,X_1}^{b,\eta ,q,Y_2,Y_1})},$$ 
the following diagrams commute. 

\begin{equation*}
\begin{CD}
Z^{\Delta }(S^1\times S^1) @>Z^{\Delta }(S)>> Z^{\Delta }(S^1\times S^1) \\ 
@V\bar{P}VV @VV\bar{P}V \\ 
Z^{\Delta }(S^1\times S^1) @>Z^{\Delta }(S')>> Z^{\Delta }(S^1\times S^1), 
\end{CD}
\quad 
\begin{CD}
Z^{\Delta }(S^1\times S^1) @>Z^{\Delta }(T^{-1})>> Z^{\Delta }(S^1\times S^1) \\ 
@V\bar{P}VV @VV\bar{P}V \\ 
Z^{\Delta }(S^1\times S^1) @>Z^{\Delta }({T'}^{-1})>> Z^{\Delta }(S^1\times S^1) 
\end{CD}
\end{equation*}

This completes the proof. 
\qed

\par \bigskip \bigskip 
\section{\kern-1em .\ \ Correspondence between Izumi's Tube Algebra and Ocneanu's Tube Algebra} 
\par 
Izumi \cite{Izumi4} translated the notion of Ocneanu's tube algebra \cite{Ocneanu1} into the language of sectors, 
and gave explicit formulas of the tube algebra operations and of the $S$- and $T$-matrices. 
He also showed that $SL(2,\mathbb{Z})$ acts on the center of the tube algebra and 
that the Verlinde identity holds. 
It is a natural question whether Izumi's tube algebra and Ocneanu's one are isomorphic as algebras, 
and whether the centers of these tube algebras are isomorphic as algebras commuting with $SL(2,\mathbb{Z})$-actions. 
In this section, we show that there exists a conjugate linear isomorphism 
between the center of Izumi's tube algebra and that of Ocneanu's one, 
which preserves products of algebras and commutes with $SL(2,\mathbb{Z})$-actions. 
\par 
Let us recall the definition of the tube algebra in sector theory, 
which was introduced by Izumi \cite{Izumi4}.  
Let $\Delta $ be a finite system of $\text{End}(M)_0$ obtained from a subfactor $N\subset M$ 
of an infinite factor $M$ with finite index and finite depth. 
The {\it tube algebra} $\text{Tube}_{I} \Delta$ introduced by Izumi is a $C^{\ast }$-algebra defined as follows. 
As a $\mathbb{C}$-vector space, $\text{Tube}_{I}\Delta$ is spanned by 
$$\{ (\xi \zeta \vert X\vert \zeta \eta )\ \vert \ \xi ,\zeta ,\eta \in \Delta ,\ X\in (\xi\zeta ,\zeta \eta )\} .$$
The product of $\text{Tube}_{I}\Delta$ is given by 
$$(\xi \zeta \vert X\vert \zeta \eta )\cdot (\xi '\zeta '\vert Y\vert \zeta '\eta ')
=\delta _{\eta ,\xi '}\sum\limits_{\nu \prec \zeta \zeta '}
\sum\limits_{A\in \mathcal{B}_{\zeta \zeta '}^{\nu }}(\xi \nu \vert A^{\ast }\zeta (Y)X\xi (A)\vert \nu \eta '),$$
where $\delta _{\eta ,\xi '}$ is Kronecker's delta and 
$\nu \prec \zeta \zeta '$ means that $\nu $ appears in the product $\zeta \zeta '$ as an irreducible component. 
The $\ast $-structure is given by 
$$(\xi \zeta \vert X\vert \zeta \eta )^{\ast }
=d(\zeta )(\eta \bar{\zeta }\vert \bar{\zeta }
(\xi (\overline{R}_{\zeta }^{\ast })X^{\ast })R_{\zeta }\vert \bar{\zeta }\xi ).$$

Izumi showed that $SL(2,\mathbb{Z})$ acts on the center $\mathcal{Z}(\text{Tube}_{I} \Delta )$ 
of $\text{Tube}_{I} \Delta $ \cite{Izumi4}. The action is given by 

$$S'_{\Delta }(\xi \eta \vert X \vert \eta \xi )
=d(\xi )(\bar{\eta }\xi \vert R_{\eta }^{\ast }\bar{\eta }
(X\xi (\bar{R}_{\eta }))\vert \xi \bar{\eta }), $$

$$T'_{\Delta }(\xi \eta \vert X \vert \eta \xi )
=\sum\limits_{\zeta \in \Delta }d(\zeta )
(\zeta \bar{\zeta }\vert R_{\zeta }\bar{R}_{\zeta }^{\ast }\vert \bar{\zeta }\zeta )
\cdot (\xi \eta \vert X \vert \eta \xi )$$
for $S'=\begin{pmatrix} 0 & -1 \\ 1 & 0
\end{pmatrix}$ and $T'=\begin{pmatrix} 1 & 1 \\ 0 & 1\end{pmatrix}$. 
In particular, the action of ${T'}^{-1}$ is given by 
$${T'}^{-1}_{\Delta }(\xi \eta \vert X \vert \eta \xi )
=(\xi p\vert X_1^{\ast }\xi (X_2)\vert p \xi )$$
for $X=X_2X_1^{\ast },\ X_1\in (p, \xi \eta ),\ X_2\in (p,\eta \xi )$ 
and $p,\xi ,\eta \in \Delta $.

\par \bigskip 
Izumi's tube algebra is isomorphic to Ocneanu's one in the following sense. 

\par \bigskip 
\begin{thm}
\label{Theorem1}  
Let $\Delta $ be a finite system of $\text{End}(M)_0$ obtained from a subfactor $N\subset M$ 
of an infinite factor $M$ with finite index and finite depth. 
Let $\varphi :\text{Tube}\Delta \longrightarrow \text{Tube}_I\Delta $ 
be the $\mathbb{C}$-linear map defined by 
$$\varphi (X_1\otimes X_2)=w^{-\frac{1}{2}}\biggl ( \frac{d(a)}{d(p)d(\rho )}\biggr ) ^{1/4}
\biggl ( \frac{d(a)}{d(p)d(\xi )}\biggr ) ^{1/4}\ (\rho a \vert X_2X_1^{\ast}\vert a \rho )
$$
for $X_1\otimes X_2\in \mathcal{H}_{\rho a}^p\otimes \mathcal{H}_{a\xi }^p$. 
Then, $\varphi $ is an algebra isomorphism. Furthermore, the restriction of $\varphi $ to the subspace 
$V^{\Delta }(S^1\times S^1)$ gives rise to an $\mathbb{C}$-linear isomorphism 
$\bar{\varphi }:Z^{\Delta }(S^1\times S^1)\longrightarrow \mathcal{Z}(\text{Tube}_I\Delta )$ 
such that the following diagram commutes for all $R\in SL(2,\mathbb{Z})$. 
\begin{equation*}
\begin{CD}
Z^{\Delta }(S^1\times S^1) @>Z^{\Delta }(R)>> Z^{\Delta }(S^1\times S^1) \\ 
@V\bar{\varphi }VV @VV\bar{\varphi }V \\ 
\mathcal{Z}(\text{Tube}_I\Delta ) @>Q(R)_{\Delta }>>  \mathcal{Z}(\text{Tube}_I\Delta )
\end{CD}
\end{equation*} 
Here, $Q:SL(2,\mathbb{Z})\longrightarrow SL(2,\mathbb{Z})$ 
is the group isomorphism satisfying $Q(S)=S'$ and $Q(T)=T'$. 
\end{thm}

\par \bigskip 
To prove the above theorem, we need the following theorem which was announced by Ocneanu \cite{Ocneanu1} 
and was proved in \cite{KSW}. 
This is an important theorem about the tube algebra. 

\par \bigskip 
\begin{thm}[\cite{Ocneanu1, KSW}]
\label{Theorem2} 
Let $\Delta $ be a finite system of $\text{End}(M)_0$ obtained from a subfactor $N\subset M$ 
of an infinite factor $M$ with finite index and finite depth. 
Then, 
$$Z^{\Delta }(S^1\times S^1)\cong \mathcal{Z}(\text{Tube}\Delta )$$
as vector spaces. 
\end{thm}

\indent 
{\bf Proof of Theorem \ref{Theorem1}.} 
First, we show that $\varphi $ is a homomorphism of algebras.  
For $X_1\otimes X_2\in \mathcal{H}_{\rho a}^p\otimes \mathcal{H}_{a\xi }^p$ and 
$Y_1\otimes Y_2\in \mathcal{H}_{\eta b}^q\otimes \mathcal{H}_{b\zeta }^q$, 
the product in $\text{Tube}\Delta $ is given by 

\begin{align*}
(X_1\otimes X_2)\cdot (Y_1\otimes Y_2)=&\delta _{\xi ,\eta }\lambda ^{-\frac{1}{2}}
\sqrt{\frac{d(a)d(b)d(r)}{d(c)d(p)d(q)d(\xi )}} \\ 
& \times \sum\limits_{\substack{Z_1\in \mathcal{B}_{\rho c}^r \\ Z_2\in \mathcal{B}_{c\zeta }^r}}Z_2^{\ast }
(\sum\limits_{A\in \mathcal{B}_{a b}^c}A^{\ast }a(Y_2Y_1^{\ast})X_2X_1^{\ast}\xi (A))Z_1\ (Z_1\otimes Z_2).
\end{align*}

\noindent 
Thus, we have 

\begin{align*}
&d(\rho )^{1/4}d(\zeta )^{1/4}\varphi ((X_1\otimes X_2)\cdot  (Y_1\otimes Y_2)) \\ 
=&\delta _{\xi ,\eta }\frac{\lambda ^{-1}}{\sqrt{d(\xi )}}\sqrt{\frac{d(a)d(b)}{d(p)d(q)}}
\sum\limits_{\substack{Z_1\in \mathcal{B}_{\rho c}^r \\ Z_2\in \mathcal{B}_{c\zeta }^r}}Z_2^{\ast }
(\sum\limits_{A\in \mathcal{B}_{a b}^c} A^{\ast }a(Y_2Y_1^{\ast})X_2X_1^{\ast}\xi (A))Z_1\ 
(\xi c\vert Z_2Z_1^{\ast } \vert c\xi ) \\ 
=&\delta _{\xi ,\eta }\frac{\lambda ^{-1}}{\sqrt{d(\xi )}}\sqrt{\frac{d(a)d(b)}{d(p)d(q)}}
\sum_{\substack{c, r\in \Delta \\ c\prec ab \\ r\prec c\xi }}
\sum\limits_{\substack{Z_1\in \mathcal{B}_{\rho c}^r \\ Z_2\in \mathcal{B}_{c\zeta }^r}}
\sum\limits_{A\in \mathcal{B}_{a b}^c}\ (\xi c\vert Z_2Z_2^{\ast }
A^{\ast }a(Y_2Y_1^{\ast})X_2X_1^{\ast}\xi (A)Z_1Z_1^{\ast } \vert c\xi )  \\ 
=&\delta _{\xi ,\eta }\frac{\lambda ^{-1}}{\sqrt{d(\xi )}}\sqrt{\frac{d(a)d(b)}{d(p)d(q)}}
\sum_{\substack{c\in \Delta \\ c\prec ab}}\sum\limits_{A\in \mathcal{B}_{a b}^c} 
(\xi c\vert A^{\ast }a(Y_2Y_1^{\ast})X_2X_1^{\ast}\xi (A) \vert c\xi ) \\ 
=&d(\rho )^{1/4}d(\zeta )^{1/4}\varphi (X_1\otimes X_2)\varphi (Y_1\otimes Y_2). 
\end{align*}

Hence, $\varphi $ preserves the products. 
It is clear that $\varphi $ is bijective. 
Thus, $\varphi $ is an algebra isomorphism. 
\par 
Let $\varphi '$ be the algebra isomorphism from the center of $\text{Tube}\Delta $ 
to the center of $\text{Tube}_I\Delta $ induced from $\varphi $. 
We regard $Z^{\Delta }(S^1\times S^1)$ as a subspace of the center $\mathcal{Z}(\text{Tube}\Delta )$, and 
denote by $\bar{\varphi }$ the restriction of $\varphi '$ to $Z^{\Delta }(S^1\times S^1)$. 
We will show that $\bar{\varphi }$ induces an $SL(2,\mathbb{Z})$-equivariant map from 
$Z^{\Delta }(S^1\times S^1)$ to $\mathcal{Z}(\text{Tube}_I\Delta )$ in the sense of the statement in the theorem. 
By Lemma \ref{Lemma2}, for $X_1\in \mathcal{B}_{\rho a}^p$ and $X_2\in \mathcal{B}_{a\rho }^p$ we have 

\begin{align*}
\tilde{S}_{\#}(X_1\otimes X_2)
&=\sum\limits_{\substack{r\in \Delta \\ Z_1\in \mathcal{B}_{\bar{a}\rho }^r\\ Z_2\in \mathcal{B}_{\rho \bar{a}}^r}}
 \kern-0.2cm \dfrac{\sqrt{d(p)d(r)}}{d(\rho )}Z_2^{\ast }\tilde{X}_2^{\ast }\bar{a}(\hat{X}_1)Z_1\ Z_1\otimes Z_2, \\ 
\tilde{T}_{\#}^{-1}(X_1\otimes X_2)
&=\sum\limits_{\substack{r\in \Delta \\ Z_1\in \mathcal{B}_{\rho p}^r\\ Z_2\in \mathcal{B}_{p\rho }^r}}
\kern-0.2cm \dfrac{\sqrt{d(a)d(r)}}{d(p)}Z_2^{\ast }X_1^{\ast }\rho(X_2)Z_1\ Z_1\otimes Z_2.
\end{align*}

It follows immediately from the second equation that 
$\varphi '\circ \tilde{T}_{\#}^{-1}={T_{\Delta }'}^{-1}\circ \varphi '$. 
Since 
 
\begin{align*}
\dfrac{1}{d(\rho )}Z_2^{\ast }\tilde{X}_2^{\ast }\bar{a}(\hat{X}_1)Z_1
&=\dfrac{1}{d(\rho )}\dfrac{d(\rho )d(a)}{d(p)}Z_2^{\ast }( \bar{a}(X_2^{\ast })R_a)^{\ast } \bar{a}(X_1^{\ast }\rho (\overline{R}_a))Z_1 \notag \\ 
&=\dfrac{d(a)}{d(p)}Z_2^{\ast }R_a^{\ast }\bar{a}(X_2)\bar{a}(X_1^{\ast }\rho (\overline{R}_a))Z_1 \\ 
&=\dfrac{d(a)}{d(p)}Z_2^{\ast }R_a^{\ast }\bar{a}(X_2X_1^{\ast }\rho (\overline{R}_a))Z_1 \\ 
&=\dfrac{d(a)}{d(p)d(\rho )}Z_2^{\ast}S_{\Delta}'(\rho a\vert X_2X_1^{\ast }\vert a\rho )Z_1, 
\end{align*}
 
\noindent 
we see also that $\varphi '\circ \tilde{S}_{\#}= S_{\Delta }'\circ \varphi '$. 
Thus, $\bar{\varphi }$ is an $SL(2,\mathbb{Z})$-equivariant homomorphism. 
\par 
Since $\dim Z^{\Delta }(S^1\times S^1)=\dim \mathcal{Z}(\text{Tube}_I\Delta )$ by Theorem \ref{Theorem2} 
and $\bar{\varphi }$ is injective, it follows that $\bar{\varphi }$ is an $SL(2,\mathbb{Z})$-equivariant isomorphism. 
This completes the proof. \qed 

\par \bigskip 
By Theorem \ref{Theorem1}, Proposition \ref{Proposition1} and Lemma \ref{Lemma3}, we have : 

\par \bigskip 
\begin{cor}
\label{Corollary}  
Let $\Delta $ be a finite system of $\text{End}(M)_0$ obtained from a subfactor $N\subset M$ 
of an infinite factor $M$ with finite index and finite depth. 
Let $\varphi :\text{Tube}\Delta \longrightarrow \text{Tube}_I\Delta $ be the algebra isomorphism 
defined as in the above theorem, and 
$P:V^{\Delta }(S^1\times S^1)\longrightarrow V^{\Delta }(S^1\times S^1)$ 
the conjugate linear isomorphism defined as in Proposition \ref{Proposition1}. 
Then, the composition $\varphi \circ P:V^{\Delta }(S^1\times S^1)\longrightarrow \text{Tube}_I\Delta $ 
induces a conjugate linear isomorphism 
$\phi :Z^{\Delta }(S^1\times S^1)\longrightarrow \mathcal{Z}(\text{Tube}_I\Delta )$, 
which preserves the products of these two algebras and commutes with the actions of $SL(2,\mathbb{Z})$. 
\end{cor}

\par \bigskip 
\begin{rem}
Izumi has introduced an inner product $\langle \ \ ,\ \ \rangle _{\text{Tube}}$ on $\text{Tube}_I\Delta $ 
as follows \cite{Izumi4}. 
$$\langle (\xi \zeta \vert X \vert \zeta \eta ),(\xi '\zeta '\vert Y \vert \zeta '\eta ') \rangle _{\text{Tube}}
=\delta _{\xi ,\xi '}\delta _{\zeta ,\zeta '}d(\xi )^2R_{\zeta }^{\ast }\bar{\zeta }(XY^{\ast })R_{\zeta }.$$ 
Let $\{ v_i\} _{i=0}^m$ be a Verlinde basis of $Z^{\Delta }(S^1\times S^1)$, 
and $\{ w_i\} _{i=1}^m$ the primitive idempotents in the fusion algebra $Z^{\Delta }(S^1\times S^1)$ 
obtained from $\{ v_i\} _{i=0}^m$ by applying the transformation $S$ (see Remarks \ref{remarks}). 
Then, the basis $\{ \frac{\sqrt{\lambda }}{d_i}\phi (w_i)\} _{i=0}^m$ is orthonormal with respect to the inner product 
$\langle \ \ ,\ \ \rangle _{\text{Tube}}$, 
where $d_i=\lambda S_{0i}$ and $Z^{\Delta }(S)v_i=\sum_j{S}_{ji}v_j$ \cite{Izumi4}.  
\par 
On the other hand, since the Verlinde basis $\{ v_i\} _{i=0}^m$ is orthonormal with respect to the inner product 
$\langle \ \ ,\ \ \rangle _{\text{TQFT}}$ based on the TQFT defined as in Section 3, 
and $Z^{\Delta }(S)$ is unitary with respect to the inner product $\langle \ \ ,\ \ \rangle _{\text{TQFT}}$, we have 

\begin{align*}
\langle \frac{\sqrt{\lambda }}{d_i}w_i,\ \frac{\sqrt{\lambda }}{d_j}w_j\rangle _{\text{TQFT}}
=&\langle \frac{1}{\sqrt{\lambda }S_{0i}}w_i,\ \frac{1}{\sqrt{\lambda }S_{0j}}w_j\rangle _{\text{TQFT}} \\ 
=&\frac{1}{\lambda }\langle Z^{\Delta }(S)^{\ast }(v_i), Z^{\Delta }(S)^{\ast }(v_j)\rangle _{\text{TQFT}}\\ 
=&\frac{1}{\lambda }\langle v_i, v_j\rangle _{\text{TQFT}} \\ 
=&\frac{1}{\lambda } \delta _{ij}. 
\end{align*}

Therefore, $\phi $ does not preserve the inner products. 
\end{rem}

\par \bigskip 
\section{\kern-1em .\ \ Calculations of Turaev-Viro-Ocneanu invariants of Basic $3$-manifolds}
\par 
Izumi explicitly gave an action 
of $SL(2, \mathbb{Z})$ on the center of the tube algebra in the language 
of sectors, and derived some formulas on Turaev-Viro-Ocneanu invariants of lens spaces 
by applying formulas in \cite{SuzukiWakui} to subfactors constructed by his method \cite{Izumi4, Izumi5}.  
In the previous section, we established a rigorous correspondence between the 
$S$- and $T$-matrices in Izumi's sector theory and the ones in Turaev-Viro-Ocneanu 
$(2+1)$-dimensional TQFT. Via this correspondence and the Dehn surgery formula 
of the Turaev-Viro-Ocneanu invariant, we can compute the Turaev-Viro-Ocneanu invariants 
of $3$-manifolds in the language of sectors. 
In this section, using techniques on sectors due to Izumi \cite{Izumi4, Izumi5}, 
we compute the Turaev-Viro-Ocneanu invariants from several subfactors for basic $3$-manifolds  
including lens spaces and Brieskorn $3$-manifolds. 
One of the most important result is that the homology $3$-sphere $S^3$ and 
the Poincar\'{e} homology $3$-sphere $\Sigma (2,3,5)$ are distinguished by the Turaev-Viro-Ocneanu invariant 
from the exotic subfactor constructed by Haagerup and Asaeda \cite{Haagerup, AsaedaHaagerup}, and 
$L(p,1)$ and $L(p,2)$ are distinguished by the Turaev-Viro-Ocneanu invariant 
from a generalized $E_6$-subfactor with $\mathbb{Z}/p\mathbb{Z}$ for $p=3,5$. 
\par \bigskip 
Let $\Delta $ be a finite system of $\text{End}(M)_0$ obtained from a subfactor $N\subset M$ 
of an infinite factor $M$ with finite index and finite depth. Since any finite dimensional $C^*$-algebra 
is semisimple, we may assume that  ${\rm Tube} \Delta =\oplus_{i=0}^r 
\text{M}_{n_i}({\mathbb C})$ as algebras, 
where $\text{M}_{n_i}(\mathbb{C})$ is the set of $n_i\times n_i$-matrices over $\mathbb{C}$. 
From each direct summand of ${\rm Tube }\Delta =\oplus_{i=0}^r 
\text{M}_{n_i}(\mathbb{C})$, we pick up a minimal projection $p_i$. 
Then, we have proved that $\{ p_i\} _{i=0}^r$ is a Verlinde basis of $Z^{\Delta }(S^1\times S^1)$ 
in the sense of Definition \ref{Definition} \cite{KSW}. Thus, we have :

\par \bigskip 
\begin{thm}[\cite{KSW}]
\label{Theorem3}
Let $\Delta $ be a finite system of $\text{End}(M)_0$ obtained from a subfactor $N\subset M$ 
of an infinite factor $M$ with finite index and finite depth. 
Then, there exists a Verlinde basis of $Z^{\Delta }(S^1\times S^1)$ in the sense of Definition \ref{Definition}. 
\end{thm}

\par \bigskip 
Let us recall the Dehn surgery formula of $Z(M)$ \cite{Wakui, KSW}. 
Let $Z$ be a $(2+1)$-dimensional TQFT, and $\{ v_i \}_{i=0}^m$ a basis of $Z(S^1 \times S^1)$.
\par 
We introduce a framed link invariant in the following way. 
Let $L=L_1\cup \cdots \cup L_r$ be a framed link with $r$-components in the 3-sphere $S^3$, 
and $h_i: D^2\times S^1\longrightarrow N(L_i)$ be the framing of $L_i$ for each $i\in \{ 1, \cdots , r\} $, 
where $N(L_i)$ denotes the tubular neighborhood of $L_i$. 
We fix an orientation for $\partial N(L_i)$ such that 
$j_i:=h_i\vert _{\partial D^2\times S^1}:S^1\times S^1\longrightarrow \partial N(L_i)$ 
is orientation preserving. 
Since the orientation for $N(L_i)$ is not compatible with the orientation for the link exterior 
$X:=\overline{S^3-N(L_1)\cup \cdots \cup N(L_r)}$, 
we can consider the cobordism parametrized boundary ${\cal W}_L:=(X; \coprod \limits_{i=1}^r j_i ,\emptyset )$ 
(See \cite{KSW} for detail). 
This cobordism induces a 
$\mathbb{C}$-linear map $Z_{{\cal W}_L}: \bigotimes \limits_{i=1}^rZ(S^1\times S^1)\longrightarrow \mathbb{C}$. 
It is easy to see that for each $i_1,\cdots ,i_r=0,1,\cdots ,m$ 
the complex number $J(L; i_1, \cdots , i_r):=Z_{{\cal W}_L}(v_{i_1}\otimes \cdots \otimes v_{i_r})$  
is a framed link invariant of $L$. 
In this setting, we proved the following proposition. 

\par\bigskip 
\begin{prop}[Dehn surgery formula \cite{KSW}]
\label{Proposition2}
Let $Z$ be a (2+1)-dimensional TQFT and $\{ v_i\} _{i=0}^m$ a basis of the fusion algebra $Z(S^1\times S^1)$ 
such that $v_0$ is the identity element in the fusion algebra. 
Let $M$ be a closed oriented 3-manifold obtained from $S^3$ by Dehn surgery along a framed link 
$L=L_1\cup \cdots \cup L_r$. 
Then, the 3-manifold invariant $Z(M)$ is given by the formula
$$Z(M)=\sum_{i_1, \cdots , i_r=0}^m S_{i_10}\cdots S_{i_r0}J(L;i_1, \cdots , i_r),$$ 
where $S_{ji} \in \mathbb{C}$, $i,j=0,1,\cdots ,m$ are defined by $Z(S) v_i=\sum_{j=0}^m S_{ji} v_j$, 
and $S:S^1\times S^1\longrightarrow S^1\times S^1$ is the orientation preserving diffeomorphism defined by 
$S(z,w)=(\bar{w},z)$ for all $(z,w)\in S^1\times S^1$. 
\end{prop}

This is a quite general formula which is a conclusion of the axioms of 
$(2+1)$-dimensional TQFT. The right-hand side of this formula looks 
very similar to the formula of the Reshetikhin-Turaev invariant of closed 
3-manifolds, although many things are missing in the above general formula, 
compared to the Reshetikhin-Turaev formula \cite{RT}. 
(See \cite{KSW, SatoWakui} for more details on the Dehn surgery formula and its applications.)

\par 
From Proposition \ref{Proposition2}, if we want to compute the Turaev-Viro-Ocneanu invariant 
of a closed 3-manifold $M$, we need to compute the $S$-matrix and the framed link 
invariants $J(L;i_1,\cdots ,i_r)$. Since we know that 
the existence of the isomorphism $Z^\Delta(S^1 \times S^1)\cong \mathcal{Z}(\text{Tube} \Delta)$ 
by Theorem \ref{Theorem3}, 
we can compute the $S$-matrix with respect 
to a Verlinde basis of $Z^{\Delta }(S^1\times S^1)$ in principle 
(see also Theorem \ref{Theorem1} and Theorem \ref{Theorem2}). 
If a 3-manifold $M$ is obtained from $S^3$ 
by Dehn surgery along an \lq\lq easy'' framed link $L$, then we can compute the 
$S$-matrix, the framed link invariants $J(L;i_1,\cdots ,i_r)$ and $Z(M)$. 
In particular, for $3$-manifolds $M$ such as lens spaces and Brieskorn 3-manifolds
 we have useful formulas for $Z(M)$ as follows. 

\begin{figure}[htbp]
\begin{center}
\setlength{\unitlength}{1cm}
\scalebox{0.9}[0.9]{\includegraphics[height=3.5cm]{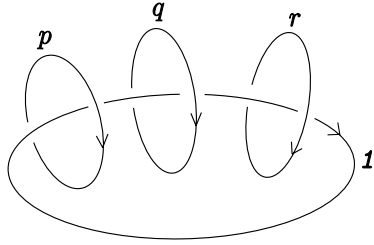}}
\caption{the Brieskorn 3-manifold $\Sigma (p,q,r)$ \label{Figure17}}
\end{center}
\end{figure}

Let $Z$ be a $(2+1)$-dimensional TQFT with Verlinde basis $\{ v_i\} _{i=0}^m$. 
For $i=0,1,\cdots ,m$, we write 
$$Z(T)v_i=t_iv_i\ \text{and} \ Z(S)v_i=\sum\limits_{j=0}^mS_{ji}v_j.$$ 
Let $p$ and $q$ be coprime positive integers. We present $p/q$ in the continued fraction 
$$\frac{p}{q}= a_1 - 
\cfrac{1}{a_2 - 
\cfrac{1}{a_3 - 
\cfrac{1}{\ddots - 
\cfrac{1}{a_n 
}}}} \ ,$$ 
where $a_1,a_2,\cdots ,a_n\geq 2$ are integers. 
Then, the lens space $L(p,q)$ is obtained by identifying two solid tori $D^2\times S^1$ 
gluing the diffeomorphism 
$f=T^{a_1}ST^{a_2}S \cdots T^{a_{n-1}}ST^{a_n}:S^1\times S^1\longrightarrow S^1\times S^1$ \cite{Rolfsen}. 
It follows from the Dehn surgery formula that 
$$Z(L(p,q))=\sum _{i_1,i_2,\cdots,i_n=0}^m
S_{i_10}t_{i_1}^{a_1}S_{i_1i_2}t_{i_2}^{a_2}S_{i_2i_3}\cdots
t_{i_{n-1}}^{a_{n-1}}S_{i_{n-1}i_n}t_{i_n}^{a_n}S_{i_n0}.$$
In particular, we have 
\par 
$\bullet$\ $Z(L(p,1))=\sum \limits_{i=0}^mt_i^pS_{i0}^2$ for any $p$. 
\par 
$\bullet$\ $Z(L(p,2))=\sum \limits_{i,j=0}^mt_i^{\frac{p+1}{2}}t_j^2S_{i0}S_{j0}S_{ij}$ for any odd integer $p$,
\par 
$\bullet$\ $Z(L(p,3))=\begin{cases} \sum \limits_{i,j,k=0}^mt_i^{\frac{p+2}{3}}t_j^2t_k^2S_{i0}S_{j0}S_{k0}S_{ij}S_{kj} & \text{if }\ p\equiv 1\ (\text{mod}3), \\  
\sum \limits_{i,j=0}^mt_i^{\frac{p+1}{3}}t_j^3S_{i0}S_{j0}S_{ij} & \text{if }\ p\equiv 2\ (\text{mod}3). \end{cases} $ 

\par \bigskip 
The Brieskorn 3-manifold 
$\Sigma (p,q,r)=\{ (u,v,w)\in \mathbb{C}^3\ \vert \ u^p+v^q+w^r=0,\ {|u|}^2+{|v|}^2+{|w|}^2=1\} $, where $p,q,r\geq 2$, 
 is obtained from $S^3$ by Dehn surgery along the framed link presented by the diagram depicted in Figure \ref{Figure17} \cite{Rolfsen}. 
 Then, it follows from the Dehn surgery formula that 
\par 
$\bullet$\ $Z(\Sigma (p,q,r))=\sum\limits_{i,j,k,l=0}^mt_i^pt_j^qt_k^rt_l\dfrac{S_{i0}S_{j0}S_{k0}S_{il}S_{jl}S_{kl}}{S_{l0}}.$
\par 
This is obtained as a special case of the following lemma. 

\vspace{1cm} \par 

\begin{figure}[hbtp]
\begin{center}
\setlength{\unitlength}{1cm}
\scalebox{0.9}[0.9]{\includegraphics[height=5cm]{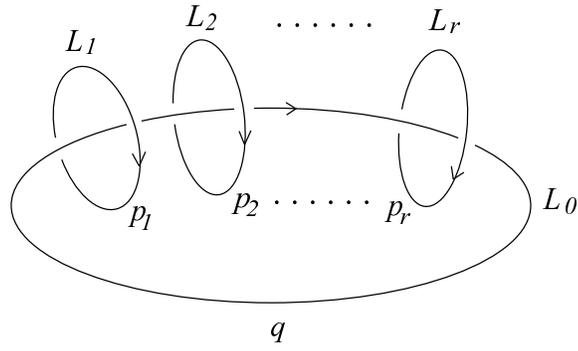}}
\caption{the framed link $L(q,p_1,...,p_r)$ \label{Figure18}}
\end{center}
\end{figure}

\par \bigskip 
\begin{lem}
\label{Lemma4}
Let $L(q,p_1,\cdots ,p_r)$ be the framed link presented by the diagram as in Figure \ref{Figure18}, and 
$M$ the $3$-manifold obtained from $S^3$ by Dehn surgery along $L(q,p_1,\cdots ,p_r)$, 
where $q,p_1,\cdots ,p_r$ stand for those integral framings.  
Then, $Z(M)$ is given by 
$$Z(M)=\sum\limits_{i_1,\cdots ,i_r,j=0}^mt_{i_1}^{p_1}\cdots t_{i_r}^{p_r}t_j^q\dfrac{S_{i_10}\cdots S_{i_r0}S_{i_1j}\cdots S_{i_rj}}{S_{j0}^{r-2}}.$$ 
\end{lem}

\par
\begin{figure}[hbtp]
\begin{center}
\setlength{\unitlength}{1cm}
\includegraphics[height=5cm]{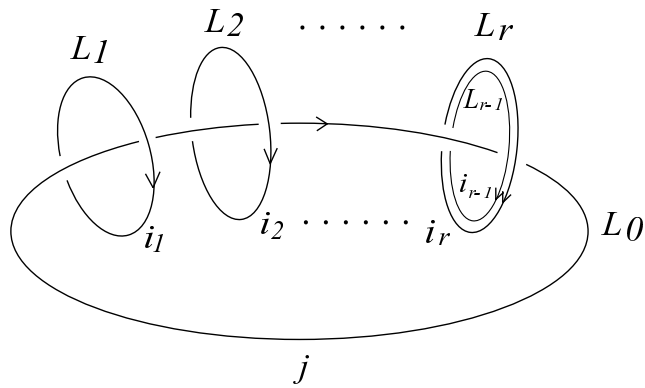}
\caption{the framed link $L(0,0,...,0)$ \label{Figure19}}
\end{center}
\end{figure}

 \medskip \indent 
{\bf Proof.} 
Let $L$ be the framed link $L(0,0,\cdots ,0)$. The framed link $L$ is isomorphic to the framed link 
presented by the diagram as in Figure \ref{Figure19}.  

By fusing the $(r-1)$-th component and $r$-th component, we obtain  

\begin{align*}
&J(L;j,i_1,\cdots ,i_{r-2},i_{r-1},i_r) \\ 
&=\sum\limits_{k=0}^mN_{i_{r-1},i_r}^kJ(L';j,i_1,\cdots ,i_{r-2},k) \\ 
&=\sum\limits_{k,l=0}^m\dfrac{S_{li_{r-1}}S_{li_r}\overline{S_{lk}}}{S_{l0}}J(L';j,i_1,\cdots ,i_{r-2},k) \\ 
&=\sum\limits_{k,l=0}^m\dfrac{S_{li_{r-1}}S_{li_r}\overline{S_{lk}}}{S_{l0}}\dfrac{S_{i_1j}\cdots S_{i_{r-2}j}S_{kj}}{S_{j0}^{r-2}} \\ 
&=\sum\limits_{l=0}^m\dfrac{S_{li_{r-1}}S_{li_r}}{S_{l0}}\dfrac{S_{i_1j}\cdots S_{i_{r-2}j}}{S_{j0}^{r-2}}\sum\limits_{k=0}^d\overline{S_{lk}}S_{kj} \qquad \qquad \notag \\ 
&=\dfrac{S_{i_1j}\cdots S_{i_{r-2}j}S_{i_{r-1}j}S_{i_rj}}{S_{j0}^{r-1}}, 
\end{align*}

\noindent 
where $L'$ is the framed link from $L$ removing the $r$-th component $L_r$. 
Hence, by induction on $r$ it can be proved that the framed link invariant $J(L;j,i_1,\cdots ,i_r)$ of 
$L=L(0,0,\cdots ,0)$ is given by 

\begin{equation}
J(L;j,i_1,\cdots ,i_r)=\dfrac{S_{i_1j}\cdots S_{i_rj}}{S_{j0}^{r-1}}.\label{eq4}
\end{equation}

By substituting the equation (\ref{eq4}) into the formula 

\begin{align*}
Z(M)
&=\sum\limits_{i_1,\cdots ,i_r=0}^mS_{i_1,0}\cdots S_{i_r,0}J(L(q,p_1,\cdots ,p_r);i_1,\cdots ,i_r) \notag \\ 
&=\sum\limits_{i_1,\cdots ,i_r=0}^mS_{i_1,0}\cdots S_{i_r,0}t_{j}^qt_{i_1}^{p_1}\cdots t_{i_r}^{p_r}J(L;i_1,\cdots ,i_r), \qquad \qquad \qquad \notag 
\end{align*}

we obtain  
$$Z(M)=\sum\limits_{i_1,\cdots ,i_r,j=0}^mt_{i_1}^{p_1}\cdots t_{i_r}^{p_r}t_j^q\dfrac{S_{i_10}\cdots S_{i_r0}S_{i_1j}\cdots S_{i_rj}}{S_{j0}^{r-2}}.$$ 
This completes the proof. \qed 

\par \bigskip 
We can compute a plenty of 
Turaev-Viro-Ocneanu invariants of some basic manifolds based on Izumi's data 
\cite{Izumi5} of the $S$- and $T$-matrices. 
For example, we have the following list of values for Turaev-Viro-Ocneanu invariants from subfactors 
by partially using the Maple software Release 5 for computations. 

\par 
\bigskip 
\noindent \underline{the $D_5^{(1)}$-subfactors \cite[Example 3.2]{Izumi2}}
\par \ 
\par 
$Z^{\text{\normalsize $D_5^{(1)}$}}(L(3,1))
=\overline{Z^{\text{\normalsize $D_5^{(1)}$}}(L(3,2))}
=\dfrac{1+2w^2}{6}$, where $w^3=1$
\par \medskip 
\par 
$Z^{\text{\normalsize $D_5^{(1)}$}}(L(5,1))
=Z^{\text{\normalsize $D_5^{(1)}$}}(L(5,2))
=Z^{\text{\normalsize $D_5^{(1)}$}}(L(7,1))
=Z^{\text{\normalsize $D_5^{(1)}$}}(L(7,2))
=\dfrac{1}{6}$

\par \bigskip 
\noindent \underline{the $E_6$-subfactor that $c_2=e^{\frac{7\pi \sqrt{-1}}{12}}$ \cite[Example 3.4]{Izumi2}}
\par \ 
\par 
$Z^{\text{\normalsize $E_6$}}(L(p,1))$
=$\dfrac{1}{12}\{ ((-1)^p+1)e^{-\frac{p\pi i}{3}}+2e^{-\frac{5p\pi i}{6}}+i^p+2(-1)^p+5\} $
\par \medskip 
$Z^{\text{\normalsize $E_6$}}(L(p,2))$
=$\dfrac{1}{4}+\dfrac{(-1)^{\frac{p+1}{2}}i}{12}-\dfrac{\sqrt{3}+i}{12}e^{-\frac{(p+1)\pi i}{6}}$
\par \medskip 
\begin{align*}
&Z^{\text{\normalsize $E_6$}}(L(p,3)) \\ 
=&
\begin{cases}
\frac{1}{24}
(9-3(-1)^{\frac{p-1}{3}}
-\sqrt{3}i(1-(-1)^{\frac{p-1}{3}})
-2\sqrt{3}i^{\frac{p-1}{3}})
&
\text{if }\ p \equiv 1 \ (\text{mod}3)\cr 
\frac{1}{24}
(9-3(-1)^{\frac{p+1}{3}}
+\sqrt{3}i(1-(-1)^{\frac{p+1}{3}})
-2\sqrt{3}i^{\frac{p+1}{3}})
&
\text{if }\ p \equiv 2 \ (\text{mod}3)
\end{cases}
\end{align*}
\par 
\bigskip 
\begin{tabular}{c|cccc}
$(p,q,r)$ &  $(2,3,5)$ & $(2,3,7)$ & $(2,5,7)$ & $(3,5,7)$ \\ 
\hline 
$Z^{\text{\normalsize $E_6$}}(\Sigma (p,q,r))$ & $\frac{2(3+\sqrt{3})+3(1-\sqrt{3})i}{12}$ & $\frac{2(3+\sqrt{3})+3(1-\sqrt{3})i}{12}$ & $\frac{-\sqrt{3}+9+6i}{12}$ & $\frac{2-\sqrt{3}i}{2}$ \\ 
\end{tabular}

\par 
\bigskip 
\begin{tabular}{c|ccccc}
$(p,q,r)$ &  $(2,3,8)$ & $(2,3,9)$ &  $(2,4,5)$ & $(2,4,6)$ &  $(2,4,7)$ \\ 
\hline 
$Z^{\text{\normalsize $E_6$}}(\Sigma (p,q,r))$ &  $\frac{(3-\sqrt{3})-(1+\sqrt{3})i}{4}$ & $\frac{3+i}{4}$ &  $\frac{2+i}{2}$ & $\frac{3+i}{2}$ & $\frac{2+i}{2}$  \\ 
\end{tabular}

\par 
\bigskip 
\begin{tabular}{c|ccccc}
$(p,q,r)$ & $(2,5,5)$ & $(2,5,6)$  & $(3,3,4)$ & $(3,3,5)$ &  $(3,3,6)$ \\ 
\hline 
$Z^{\text{\normalsize $E_6$}}(\Sigma (p,q,r))$ & $\frac{3+i}{4}$  & $\frac{3-\sqrt{3}i}{3}$ & $\frac{\sqrt{3}(1-i)}{4}$ & $\frac{-\sqrt{3}i}{2}$& $\frac{3-\sqrt{3}}{4}i$ \\ 
\end{tabular}

\par 
\bigskip 
\begin{tabular}{c|ccc}
$(p,q,r)$ & $(3,4,5)$ & $(4,4,4)$  & $(3,4,4)$ \\ 
\hline 
$Z^{\text{\normalsize $E_6$}}(\Sigma (p,q,r))$ &  $\frac{2(\sqrt{3}+3)+3(1-\sqrt{3})i}{12}$ & $\frac{3i-\sqrt{3}}{3}$  &  $\frac{3(1+\sqrt{3})+(3-\sqrt{3})i}{6}$ \\ 
\end{tabular}

\par \bigskip 
\noindent \underline{a generalized $E_6$-subfactor with $G=\mathbb{Z}/3\mathbb{Z}$ \cite[Example A-1]{Izumi5}}
\par \ 
\par 
$Z^{\text{\normalsize $E_6,\mathbb{Z}/3\mathbb{Z}$}}(L(3,1))
=\overline{Z^{\text{\normalsize $E_6,\mathbb{Z}/3\mathbb{Z}$}}(L(3,2))}
=\frac{(7-\sqrt{7}i)(\sqrt{21}-1)}{70}$
\par \medskip 
$Z^{\text{\normalsize $E_6,\mathbb{Z}/3\mathbb{Z}$}}(L(5,1))
=Z^{\text{\normalsize $E_6,\mathbb{Z}/3\mathbb{Z}$}}(L(5,2))
=-\frac{2}{15}+\frac{4\sqrt{21}}{105}$
\par \medskip 
$Z^{\text{\normalsize $E_6,\mathbb{Z}/3\mathbb{Z}$}}(L(7,1))
=Z^{\text{\normalsize $E_6,\mathbb{Z}/3\mathbb{Z}$}}(L(7,2))
=\frac{(1+\sqrt{3}i)(\sqrt{21}-1)}{30}$

\par 
\par \bigskip 
\noindent \underline{a generalized $E_6$-subfactor with $G=\mathbb{Z}/4\mathbb{Z}$ \cite[ Example A-2]{Izumi5}}
\par \ 
\par 
$Z^{\text{\normalsize $E_6,\mathbb{Z}/4\mathbb{Z}$}}(L(3,1))
=Z^{\text{\normalsize $E_6,\mathbb{Z}/4\mathbb{Z}$}}(L(3,2))
=\frac{2+\sqrt{2}}{16}$
\par \medskip 
$Z^{\text{\normalsize $E_6,\mathbb{Z}/4\mathbb{Z}$}}(L(5,1))
=Z^{\text{\normalsize $E_6,\mathbb{Z}/4\mathbb{Z}$}}(L(5,2))
=\frac{2+\sqrt{2}}{16}$
\par \medskip 
$Z^{\text{\normalsize $E_6,\mathbb{Z}/4\mathbb{Z}$}}(L(7,1))
=Z^{\text{\normalsize $E_6,\mathbb{Z}/4\mathbb{Z}$}}(L(7,2))
=\frac{2-\sqrt{2}}{16}$

\bigskip 
\noindent \underline{a generalized $E_6$-subfactor with $G=\mathbb{Z}/2\mathbb{Z}\times \mathbb{Z}/2\mathbb{Z}$ 
\cite[Example A-3]{Izumi5}}
\par \ 
\par 
$Z^{\text{\normalsize $E_6,\mathbb{Z}/2\mathbb{Z}\times \mathbb{Z}/2\mathbb{Z}$}}(L(3,1))
=Z^{\text{\normalsize $E_6,\mathbb{Z}/2\mathbb{Z}\times \mathbb{Z}/2\mathbb{Z}$}}(L(3,2))
=\frac{2+\sqrt{2}}{16}$
\par \medskip 
$Z^{\text{\normalsize $E_6,\mathbb{Z}/2\mathbb{Z}\times \mathbb{Z}/2\mathbb{Z}$}}(L(5,1))
=Z^{\text{\normalsize $E_6,\mathbb{Z}/2\mathbb{Z}\times \mathbb{Z}/2\mathbb{Z}$}}(L(5,2))
=\frac{2+\sqrt{2}}{16}$
\par \medskip 
$Z^{\text{\normalsize $E_6,\mathbb{Z}/2\mathbb{Z}\times \mathbb{Z}/2\mathbb{Z}$}}(L(7,1))
=Z^{\text{\normalsize $E_6,\mathbb{Z}/2\mathbb{Z}\times \mathbb{Z}/2\mathbb{Z}$}}(L(7,2))
=\frac{2-\sqrt{2}}{16}$

\par \bigskip 
\noindent \underline{a generalized $E_6$-subfactor with $G=\mathbb{Z}/5\mathbb{Z}$ \cite[Example A-4]{Izumi5}}
\par \ 
\par 
$Z^{\text{\normalsize $E_6,\mathbb{Z}/5\mathbb{Z}$}}(L(3,1))
=Z^{\text{\normalsize $E_6,\mathbb{Z}/5\mathbb{Z}$}}(L(3,2))
=\frac{\sqrt{5}-7}{55}$
\par \medskip 
$Z^{\text{\normalsize $E_6,\mathbb{Z}/5\mathbb{Z}$}}(L(5,1))
=\dfrac{1}{2}Z^{\text{\normalsize $E_6,\mathbb{Z}/5\mathbb{Z}$}}(L(5,2))
=\frac{1+3\sqrt{5}}{33}$
\par \medskip 
$Z^{\text{\normalsize $E_6,\mathbb{Z}/5\mathbb{Z}$}}(L(7,1))
=Z^{\text{\normalsize $E_6,\mathbb{Z}/5\mathbb{Z}$}}(L(7,2))
=\frac{3}{55}+\frac{\sqrt{5}}{33}$

\par \bigskip 
\noindent \underline{the Haagerup subfactor of Jones index $\frac{5+\sqrt{13}}{2}$ \cite[Appendix C]{Izumi5}}
\par \ 
\par 
$Z^{Haagerup}(L(3,1))=Z^{Haagerup}(L(3,2))=\frac{13-\sqrt{13}}{26}$
\par \medskip 
$Z^{Haagerup}(L(5,1))=Z^{Haagerup}(L(5,2))=\frac{13+3\sqrt{13}}{78}$
\par \medskip 
$Z^{Haagerup}(L(7,1))=Z^{Haagerup}(L(7,2))=\frac{13+3\sqrt{13}}{78}$
\par \medskip 
$Z^{Haagerup}(\Sigma (2,3,5))=-\frac{\sqrt{13}}{26}+\frac{7}{6}=Z^{Haagerup}(S^3)+1$

\par \bigskip 
From the above results of computations, we expect that the following conjecture will hold true. 

\par 
\begin{conj} If there exists a generalized $E_6$-subfactor with group symmetry $G=\mathbb{Z}/7\mathbb{Z}$, 
then the lens spaces $L(7,1)$ and $L(7,2)$ will be distinguished by the Turaev-Viro-Ocneanu invariant from the subfactor. 
\end{conj}

\end{document}